\documentclass{aastex}
\usepackage{spr-astr-addons}
\usepackage{url}

\RequirePackage{color}

\newcommand{\emaila}{jbg84@xanum.uam.mx}

\newtheorem{teo}{Theorem}
\newtheorem{dfn}{Definition}
\begin{document}

\title{Periodic orbits in the restricted four--body problem \\ with two equal masses}
\shorttitle{Periodic orbits in  restricted four--body problem}
\shortauthors{Burgos and Delgado}

\author{Jaime Burgos--Garc\'ia}\and
\author{Joaqu\'in Delgado}
\affil{Departamento de Matem\'aticas UAM--Iztapalapa.
Av. San Rafael Atlixco 186, Col. Vicentina, C.P.
09340, M\'exico, D.F.}
\email{\emaila}
\affil{Accepted for publication in Astrophysics \& Space Science}2




\begin{abstract}
 The restricted (equilateral) four-body problem consists of three bodies of masses $m_{1}$, $m_{2}$ and $m_{3}$ (called primaries) lying in a Lagrangian configuration of the three-body problem i.e., they remain fixed at the apices of an equilateral triangle in a rotating coordinate system. A massless fourth body moves under the
Newtonian gravitation law due to the three primaries; as in the restricted three-body problem (R3BP), the fourth mass does not affect the motion of the three primaries. In this paper we explore symmetric periodic orbits of the restricted four-body problem (R4BP) for the case of two equal masses where they satisfy approximately the Routh's critical value. We will classify them in nine families of periodic orbits. We offer an exhaustive study of each family and the stability of each of them.
\end{abstract}



\section{Introduction}
Few bodies problems have been studied for long time in celestial
mechanics,  either as simplified models of more complex planetary
systems or as benchmark models where new mathematical theories can
be tested. The three--body problem has been source of inspiration
and study in Celestial Mechanics since Newton and Euler. In recent
years it has been discovered multiple stellar systems such as double
stars and triple systems. The restricted three body problem (R3BP) has
demonstrated to be a good model of several systems in our solar
system such as the Sun-Jupiter-Asteroid system, and with less
accuracy the Sun-Earth-Moon system.  In analogy with the R3BP, in
this paper we  study a restricted problem of four bodies consisting
of three primaries moving in circular orbits keeping an equilateral
triangle configuration and a massless particle moving under the
gravitational attraction of the primaries. Here we focus on the
study of families of periodic orbits. We refer to this as the
restricted four body problem (R4BP). There exist some preliminary
studies of this problem in different versions, \cite{Simo}, \cite{Lea}, \cite{Ped} and \cite{PapaII} studied the equilibrium points
and their stability of this problem. Other authors have studied the
case where the primaries form a collinear configuration. At the time
of writing this paper we became aware of the paper \cite{PapaI}, where they performed a numerical
study similar to ours for two cases depending on the masses of  the
primaries: (a) three equal masses and (b) two equal masses. It is
the second case that our work is related to \cite{PapaI}, although
we use a slightly different value of the mass parameter. The reason
is the  same  as the cited authors, of having the primaries moving in
linearly stable circular orbits for a value of the mass parameter less but approximately equal to
Routh's critical value. By historical and theoretical aspects, we use the same letters used in the Copenhagen category of the R3BP to denote the families of periodic orbits, see \cite{Sz}. The families $g$, $f$, $a$, $m$, $r_{2}$, $g_{4}$, $g_{6}$ are similar to those families denoted by the same letter in the R3BP, i.e., the family of direct periodic orbits around the mass $m_{1}$ of this paper is denoted by the letter $g$ as it was done in the Copenhagen category for each family, but the families $j$ and $j_{2}$ are exclusive of this problem because they do not have similar families in the R3BP.
 Our results confirm and extend four families of periodic orbits found in
\cite{PapaI}, such families are $f_{5}$, $f_{7}$, $f_{2}$, $f_{3}$. These families correspond respectively to the first phases (defined in 3.1) of the families $a$, $j$, $f$, $g$ of this paper, however we present 5 new families of periodic orbits.
We used systematically regularization of binary
collisions of the infinitesimal with any of the primaries by a
method similar to Birkhoff's which permit us to continue some of the
families beyond double collisions. In this way we can show that such
continued families end up in a homoclinic connection. This last
phenomenon can be dynamically explained by the so called blue sky
catastrophe and a rigorous justification will appear elsewhere.

We recall that Routh's criterion for linear stability  states that
\begin{displaymath}
\frac{m_{1}m_{2}+m_{2}m_{3}+m_{3}m_{1}}{m_{1}+m_{2}+m_{3}}<\frac{1}{27}.
\end{displaymath}
When the three masses are such that $m_{2}=m_{3}:=\mu$ and
$m_{1}+m_{2}+m_{3}=1$, the inequality is satisfied in the interval
$\mu\in[0,0.019063652805978857\ldots)$, so in our case we take the masses equal to
$m_{2}=m_{3}=0.0190636$ and $m_{1}=0.9809364$.

\section{Equations of Motion}
\begin{figure}[!hbp]
\centering
\includegraphics[width=0.45\textwidth]{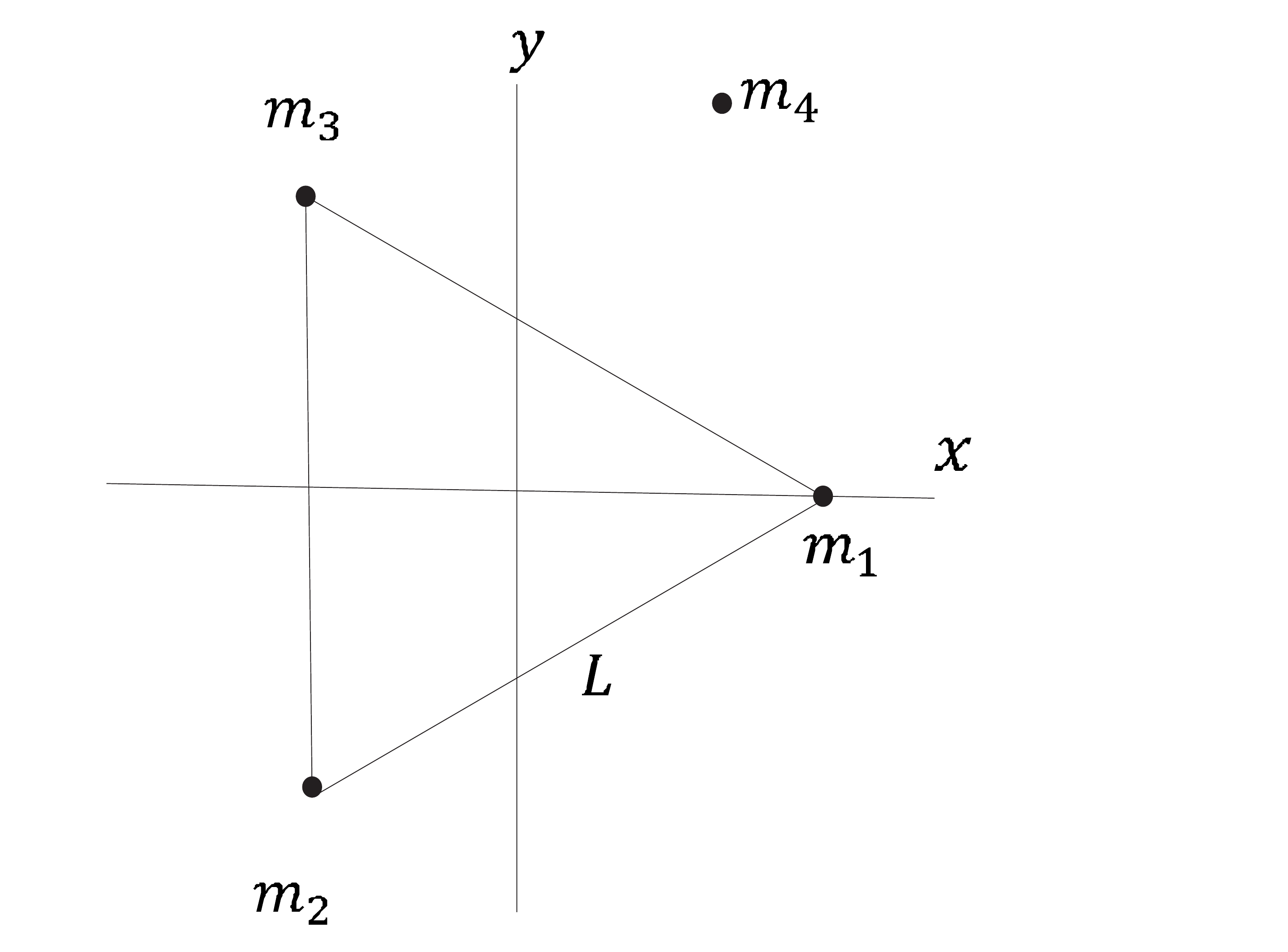}
\caption{The restricted four-body problem in a synodic system\label{triangle}}
\end{figure}

Consider three point masses, called $\textit{primaries}$, moving in
circular periodic orbits around their center of mass under their
mutual Newtonian gravitational attraction,  forming an equilateral
triangle configuration. A third massless particle moving in the same
plane is acted upon the attraction of the primaries. The equations
of motion of the massless particle referred to a synodic frame with
the same origin, where the primaries remain fixed, are:
\begin{eqnarray}
\bar{x}''-2n\bar{y}'-n^2\bar{x}&=&-k^2\sum_{i=1}^{3}m_{i}\frac{(\bar{x}-\bar{x_{i}})}{\rho_{i}^3}\nonumber\\
\bar{y}''+2n\bar{x}'-n^2\bar{y}&=&-k^2\sum_{i=1}^{3}m_{i}\frac{(\bar{y}-\bar{y_{i}})}{\rho_{i}^3} \label{sistemaconunidades}
\end{eqnarray}
where $k^2$ is the gravitational constant, $n$ is the mean motion,
$\rho_{i}^{2}=(\bar{x}-\bar{x}_{i})^2+(\bar{y}-\bar{y}_{i})^2$ is the distance of the massless particle to the primaries,
 $\bar{x}_{i}$, $\bar{y}_{i}$
are the vertices of equilateral triangle  formed by the primaries,
and ($'$)  denotes derivative with respect to time $t^{*}$. We
choose the orientation of the triangle of masses such that $m_1$
lies along the positive $x$--axis and $m_2$, $m_3$ are located
symmetrically with respect to the same axis, see
Figure~\ref{triangle}.

The equations of motion can be recast in dimensionless form as follows:
Let $L$ denote the length of triangle formed by the primaries,  $x=\bar{x}/L$,
$y=\bar{y}/L$, $x_i=\bar{x}_i/L$, $y_i=\bar{y}_i/L$, for $i=1,2,3$;
$M=m_{1}+m_{2}+m_{3}$ the total mass, and  $t=nt^{*}$. Then the equations (\ref{sistemaconunidades}) become
\begin{eqnarray}
\ddot{x}-2\dot{y}-x &=&-\sum_{i=1}^{3}\mu_{i}\frac{(x-x_{i})}{r_{i}^3}\nonumber\\
\ddot{y}+2\dot{x}-y &=&-\sum_{i=1}^{3}\mu_{i}\frac{(y-y_{i})}{r_{i}^3} \label{sistemasinunidades}
\end{eqnarray}
where we have used Kepler's third law: $k^2M=n^2L^3$, ($\dot{}$) represents derivatives with respect to the dimensionless
time $t$ and $r_{i}^2=(x-x_{i})^2+(y-y_{i})^2$.

 The system
(\ref{sistemasinunidades}) will be defined if we know the vertices
of triangle for each value of the masses. In this paper we suppose
$\mu:=m_{2}=m_{3}$ then $\mu_{1}=1-2\mu$. It is  not difficult to show
that the vertices of triangle are given by   $x_{1}=\sqrt{3}\mu$, $y_{1}=0$,
$x_{2}=-\frac{\sqrt{3}(1-2\mu)}{2}$, $y_{2}=-\frac{1}{2}$,
$x_{3}=-\frac{\sqrt{3}(1-2\mu)}{2}$, $y_{3}=\frac{1}{2}$. The system
(\ref{sistemasinunidades}) can be written succinctly as

\begin{eqnarray}
\ddot{x}-2\dot{y}&=&\Omega_{x} \label{sistemastandar}\\
\ddot{y}+2\dot{x}&=&\Omega_{y}
\end{eqnarray}
where
\begin{displaymath}
\label{omega}\Omega(x,y,\mu):=\frac{1}{2}(x^2+y^2)+\sum_{i=1}^{3}\frac{\mu_{i}}{r_{i}}
\end{displaymath}
is the effective potential function.

 There are three limiting cases:
\begin{enumerate}
\item If $\mu=0$, we obtain the rotating Kepler's problem, with $m_{1}=1$ at the origin of coordinates.
\item If $\mu=1/2$, we obtain the circular restricted three body problem, with two equal masses $m_{2}=m_{3}=1/2$.
\item If $\mu=1/3$, we obtain the symmetric case with three masses equal to $1/3$.
\end{enumerate}

It will be useful to write the system (\ref{sistemastandar}) using
complex notation. Let $z=x+\textit{i}y$, then
\begin{equation}
\label{sistemacomplejo}\ddot{z}+2\textit{i}\dot{z}=2\frac{\partial\Omega}{\partial\bar{z}}
\end{equation}
with
\begin{displaymath}
\Omega(z,\bar{z},\mu)=\frac{1}{2}\vert
z\vert^2+U(z,\bar{z},\mu)
\end{displaymath}
where the gravitational potential is
\begin{displaymath}
U(z,\bar{z},\mu)=\sum_{i=1}^{3}\frac{\mu_{i}}{\vert z-z_{i}\vert}
\end{displaymath}
and $r_{i}=\vert z-z_{i}\vert$, $i=1,2,3$ are the distances to the
primaries. System (\ref{sistemacomplejo}) has the Jacobian first
integral
\begin{displaymath}
2\Omega(z,\bar{z},\mu)-\vert\dot{z}\vert^{2}=C.
\end{displaymath}

If we define $P=p_{x}+\textit{i}p_{y}$, the conjugate momenta of
$z$, then system (\ref{sistemastandar}) can be recast as  a
Hamiltonian system  with Hamiltonian
\begin{eqnarray}
H&=&\frac{1}{2}\vert P\vert^2+Im(z\overline{P})-U(z,\bar{z},\mu)\nonumber\\
&=&\frac{1}{2}(p^{2}_{x}+p^{2}_{y})+(yp_{x}-xp_{y})-U(x,y,\mu).\label{hamiltoniano}
\end{eqnarray}
The relationship with the Jacobian integral is $H=-C/2$.
The phase space of (\ref{hamiltoniano}) is defined as
\begin{displaymath}
\Delta=\{(z,P)\in\mathbb{C}\times\mathbb{C}\vert z\ne z_{i},
i=1,2,3\},
\end{displaymath}
with collisions occurring at $z=z_{i}$, $i=1,2,3$.

\begin{figure}[h]
\includegraphics[width=0.2\textwidth]{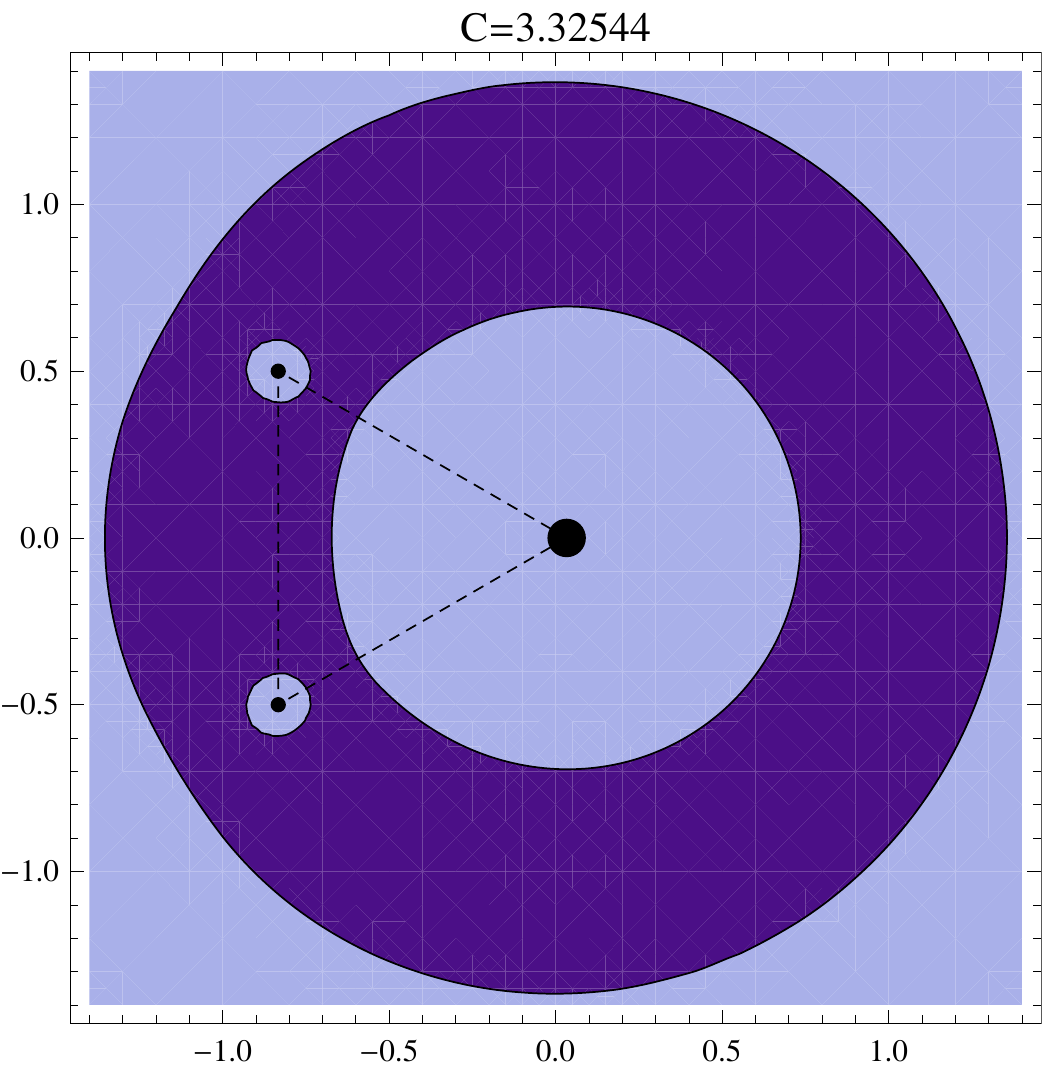}
\includegraphics[width=0.2\textwidth]{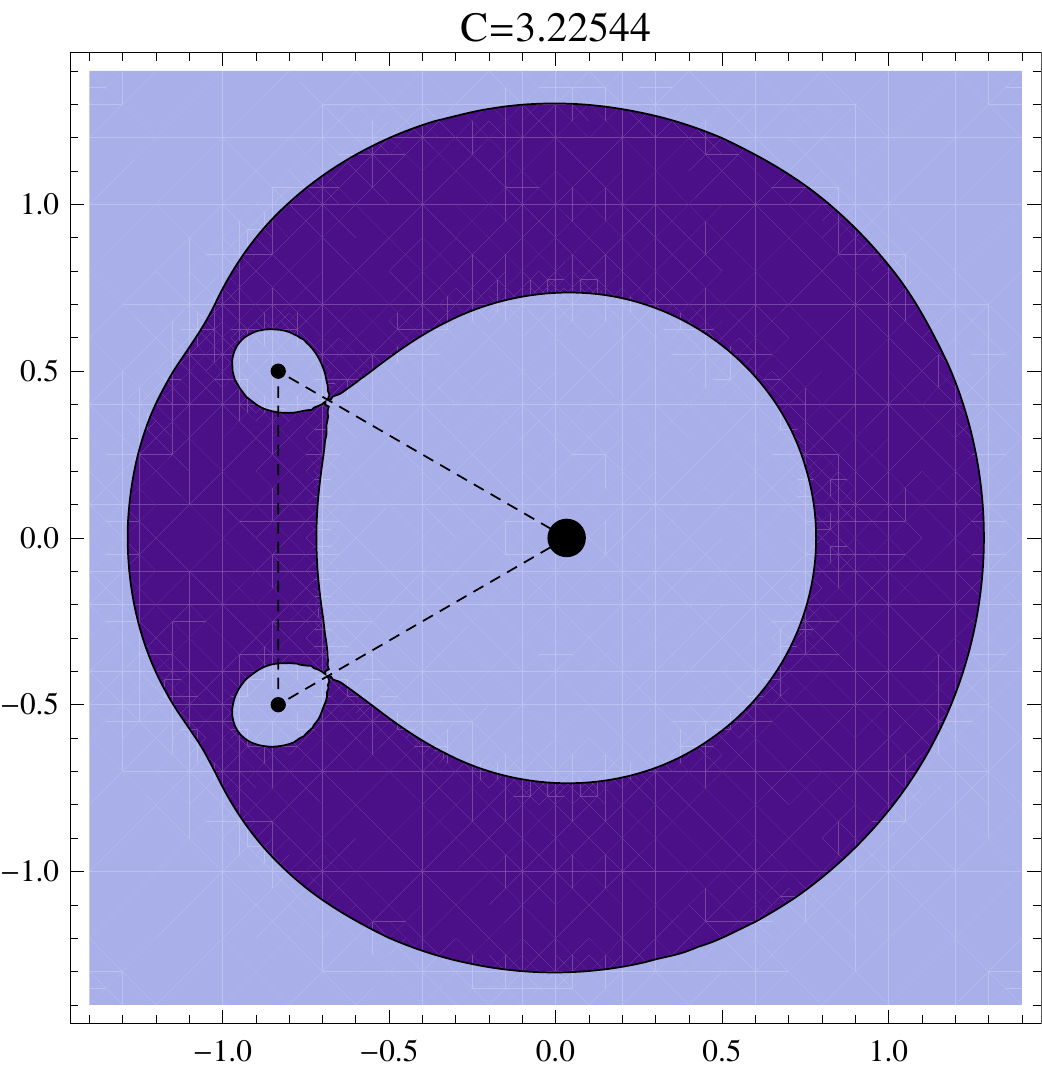}\\
\includegraphics[width=0.2\textwidth]{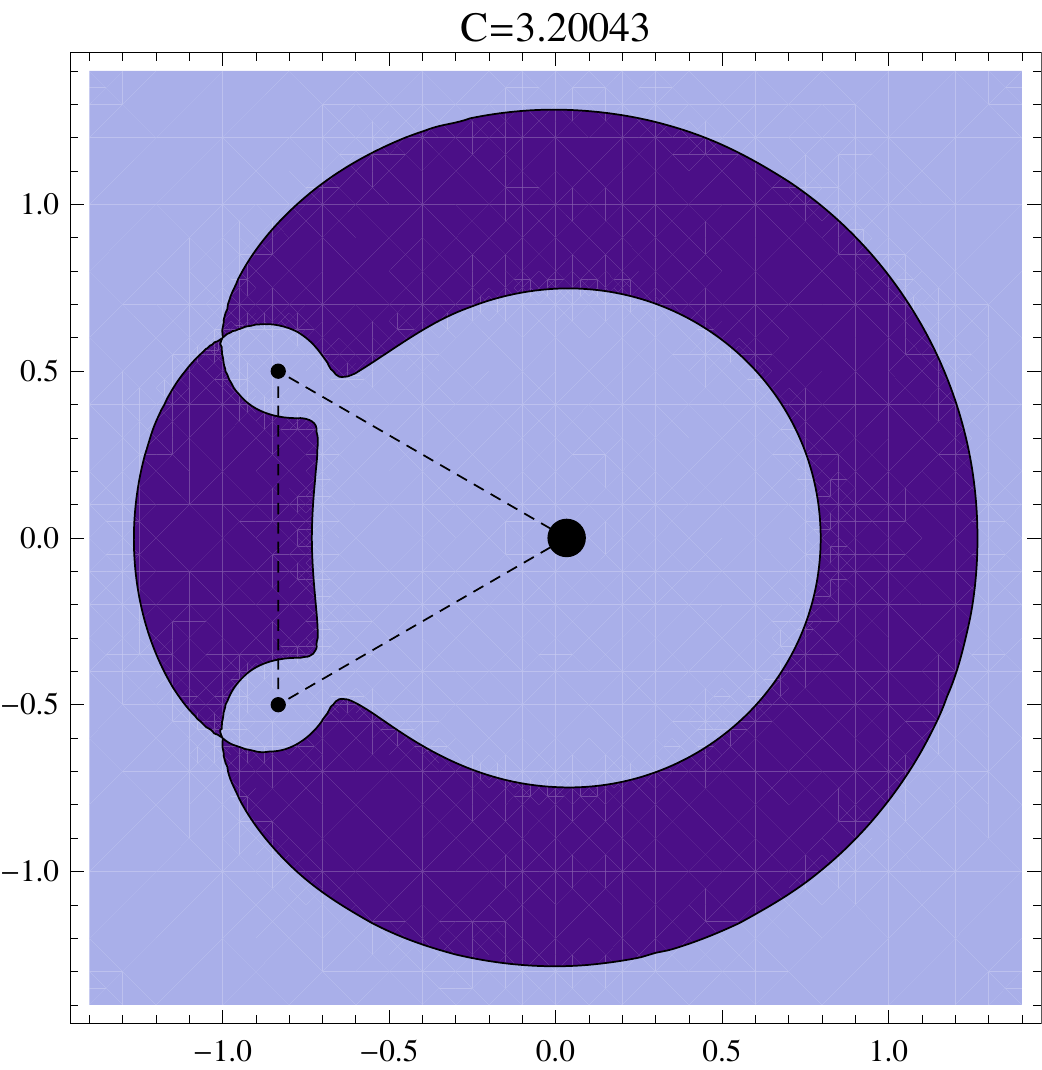}
\includegraphics[width=0.2\textwidth]{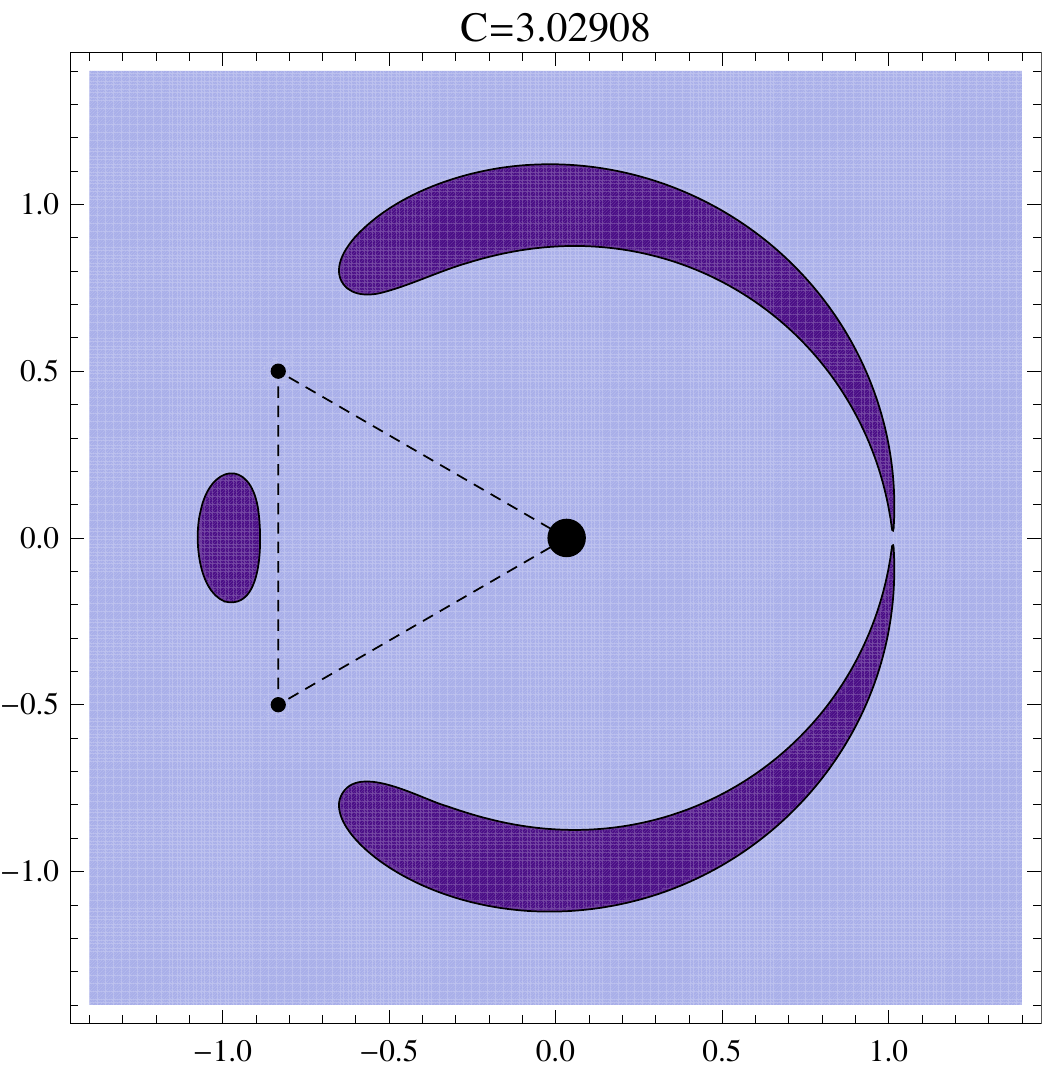}\\
\includegraphics[width=0.2\textwidth]{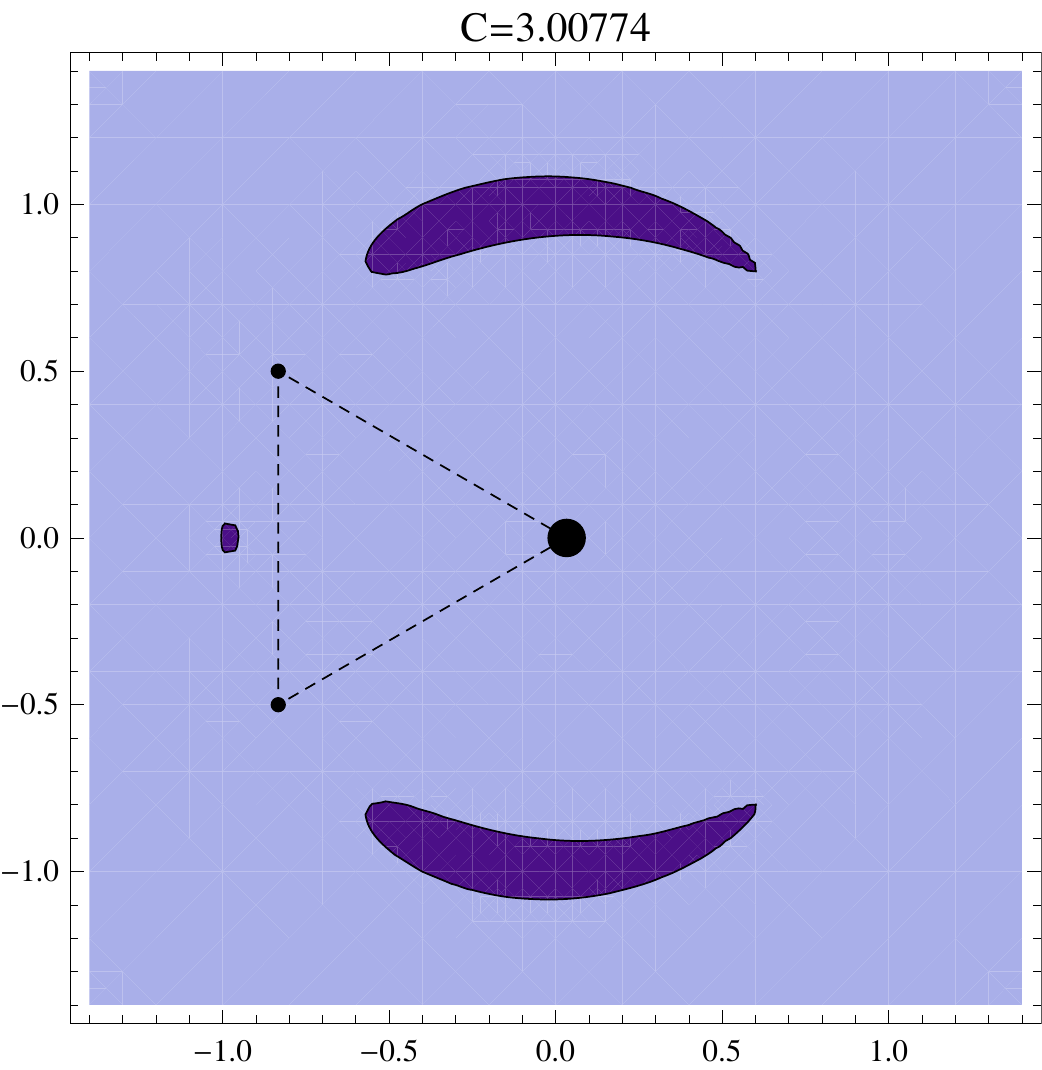}
\includegraphics[width=0.2\textwidth]{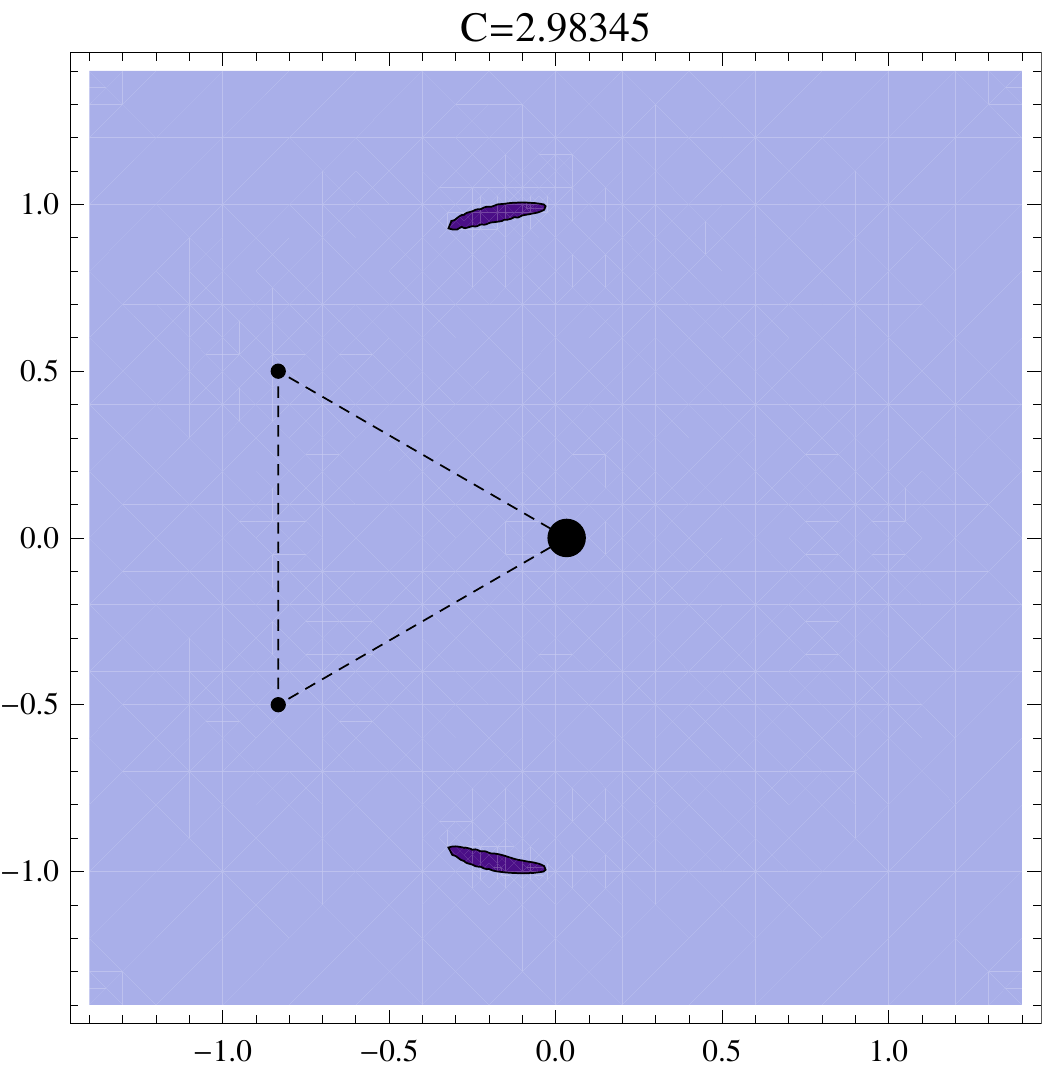}
      \caption{Hill's regions for a large value of the Jacobian constant (top-left). Hill's regions for critical values of the Jacobian constant (top-right and second row). The Hill's regions of the last row correspond to a slightly larger value than the critical one for illustrations purposes. }\label{others}
\end{figure}

There exist five equilibrium points for all values of the masses of the primaries in the R3BP;
in the R4BP, the number of equilibrium
points depends on the particular values of the masses. For the value of the mass parameter we are using throughout this paper of $\mu=0.0190636$, the Hill's regions are shown in Figure~\ref{others}. For large values of $C$ the Hill's regions consist of small disks around the primaries together with an unbounded component having as boundary a closed curve around the primaries. As the Jacobian constant decreases, the evolution of the Hill's region is shown in Figure~\ref{others}. The smaller value of $C$ is just above the critical value where the Hill's region is the whole plain minus the positions of the primaries.

\begin{figure}[h]
 \centering
\includegraphics[width=0.35\textwidth]{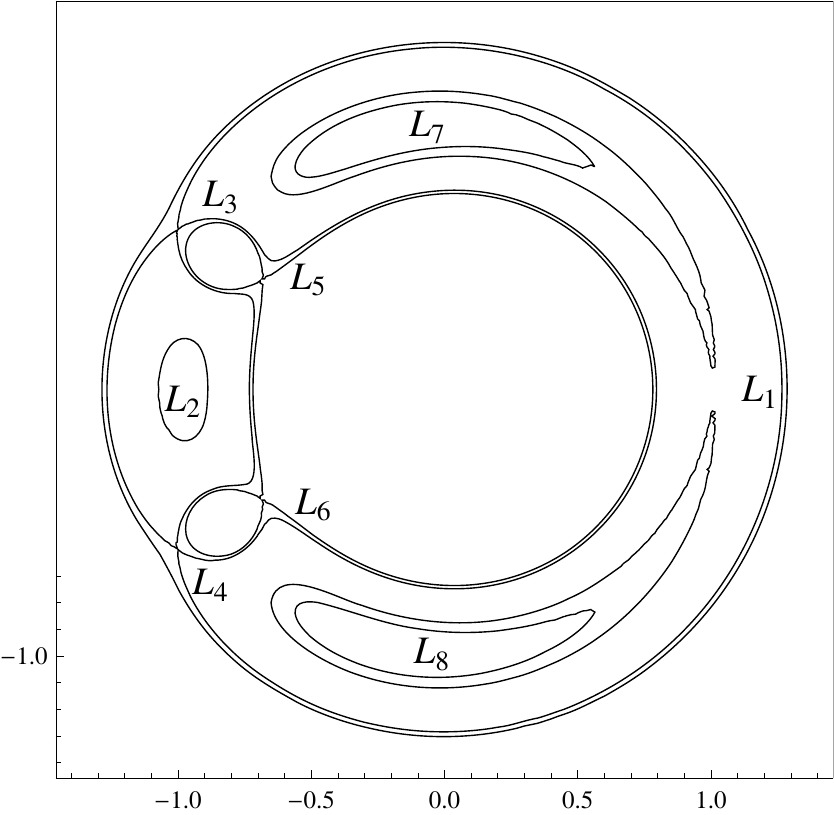}
\caption{The eight equilibrium points (green) for two equal masses.\label{eqpoints}}
\end{figure}

 A complete discussion of the equilibrium points and
bifurcations can be found in \cite{Del}, \cite{MeyerCC}, \cite{Lea},  \cite{PapaII}, \cite{Simo}.
In our particular problem we have 2 collinear and 6 non-collinear equilibrium points.  We use the notation shown in Figure~\ref{eqpoints} for the eight critical points. All  of them are unstable except the non-collinear  $L_{7}$ and $L_{8}$.

\section{Symmetric periodic orbits}
In what follows we consider symmetric periodic orbits, i,e.,
periodic orbits symmetric with respect to the synodical $x$-axis,
namely orbits which are invariant  under  the symmetry
$(z,t)\rightarrow(\bar{z},-t)$. Thus  a symmetric periodic orbit is
defined by   two successive perpendicular crossings with the
$x$-axis.  We use as staring point of the continuation either nearly
Keplerian circular orbits for large values of the Jacobian constant
or small Liapunov's orbits emerging form  $L_{1}$. It was shown
in \cite{Lea} that the equilibrium point
$L_{1}$ is unstable with eigenvalues $\pm\textit{i}\omega$,
$\pm\lambda$ where $\omega$ and $\lambda$ are real numbers, so we
can use the Liapunov's Center Theorem to find periodic orbits around
the equilibrium point $L_{1}$ and continue them.

In the R3BP the equilibrium points $L_{4}$ and $L_{5}$ are limits of
families of periodic orbits of the family $g$ of the Copenhagen
category. This phenomenon  is known as the ``blue sky catastrophe"
termination principle. In the papers by \cite{Buf}, \cite{MeyMcSw} and \cite{Simo} a
complete discussion of this principle in the R3BP is given. We state
the theorem behind this phenomena as stated in \cite{Henr}:

\begin{teo}  \label{blue sky} Let us consider a non-degenerate homoclinic orbit to an equilibrium with eigenvalues $\pm\alpha\pm i\omega$ with $\alpha$ and $\omega$ reals and strictly positive, of a real analytic Hamiltonian system with two degrees of freedom. Close to this orbit there exists an analytical family of periodic orbits with the following properties
\begin{enumerate}

\item The family can be parametrized by a parameter $\epsilon$ in the interval $0<\epsilon<\delta$ with $\delta$ small enough. Let us write it as $x(\epsilon,t)=x(\epsilon,t+T(\epsilon))$.
\item For $\epsilon=0$ we have the homoclinic orbit.
\item The period $T(\epsilon)$ increases without bound when $\epsilon$ goes to zero.
\item The characteristic exponents of the family change from the stable type to unstable type an vice-versa infinitely many times as $\epsilon$ goes to zero.
\end{enumerate}

\end{teo}
In a future work we will discuss the theoretical aspects of this termination principle in the R4BP.

In the following section we will classify a large number of periodic orbits in sets called families. We will take the Henon's definition for a family of periodic orbits to make such classification (see \cite{HenGen}):

\begin{dfn} A family of periodic orbits is a set of symmetric periodic orbits for which the initial parameter $x_{0}=x(0)$ and the period in family can be considered as two continuous functions of one single parameter $\alpha$.
\end{dfn}

In general we will consider the Str\"{o}mgren's termination principle to decide when a family ends.
\begin{dfn} Suppose we have obtained a finite section of a family of periodic orbits in an interval of the parameter $[\alpha_{1},\alpha_{2}]$ for which the family is followed and we want to extend it, then
\begin{enumerate}
\item The family remains in itself, i.e., the characteristic curve is a closed curve. We call this family a closed family.
\item For $\alpha\leq\alpha_{1}$ and $\alpha\geq\alpha_{2}$ the family has a \textit{natural termination} for which one of the following amounts grow without limit
\begin{itemize}
\item The dimension $D$ of the orbit, defined as the maximum distance to the origin.
\item The parameter $\alpha$.
\item The period of the orbit.
\end{itemize}
\end{enumerate}
The second case is called an open family.
\end{dfn}
 Note that the principle of termination of a family of periodic orbits mentioned in \ref{blue sky} is a particular case of the above definition because the period of the orbit increases without bound.
\subsection{The search for periodic orbits}
The periodic orbits were calculated in double precision with a multi-step Adams-Bashforth integrator of variable order for more accuracy. New transformations were needed to regularize different kind of collisions appearing in the families of periodic orbits of this problem. The families have been identified by letters as in the Copenhagen category with or without subscripts, the subscripts meaning the number of loops of the orbit around the primary under consideration. We use the classic $(x,C)$ plane of characteristic curves to represent the families of periodic orbits, in addition we use the $(a,C)$ plane to show the evolution of the stability of the families, here $a$ denotes the stability index (see \cite{HenII} for details), we have stability in the linear sense when $\vert a\vert<1$ and instability in other case.

The families have been separated in phases as in the R3BP (see \cite{Sz}), the colors in the characteristic curves of the nine families represent the different phases (and orbits near to collision) of each family, representative orbits are shown to illustrate each phase. Some orbits shown by \cite{Brou} and \cite{PapaI} can be compared with ours.
\section{Classification of families of periodic orbits}
\begin{figure}[!h]
\begin{center}
\includegraphics[width=0.45\textwidth,height=0.45\textwidth]{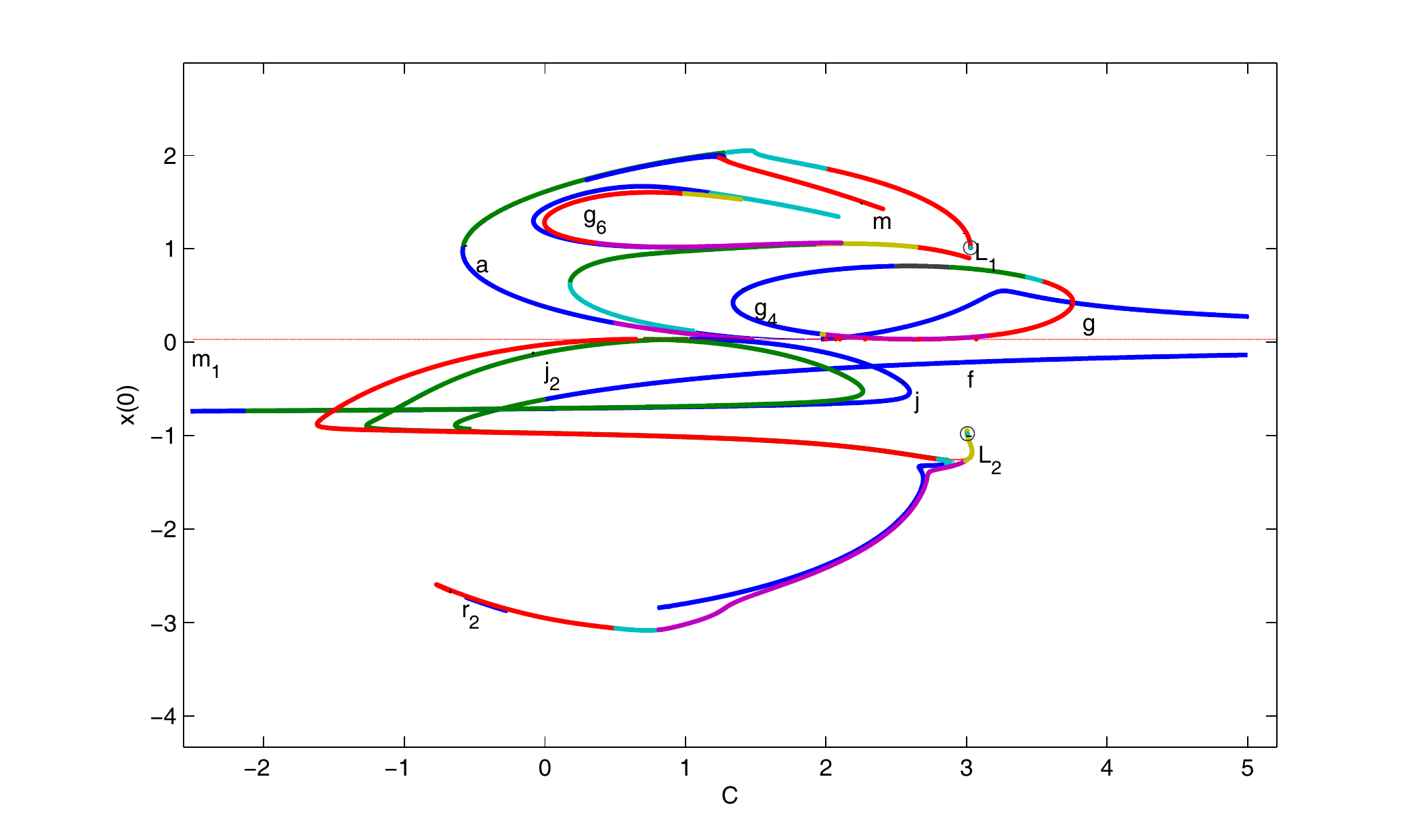}
\end{center}
\caption{The characteristic curves of nine families of periodic orbits for the restricted four-body problem with two equal masses.\label{characteristics}}
\end{figure}
\subsection{The family $g$ of direct periodic orbits around $m_{1}$}
\begin{figure}[!h]
\begin{center}
\includegraphics[width=0.45\textwidth,height=0.45\textwidth]{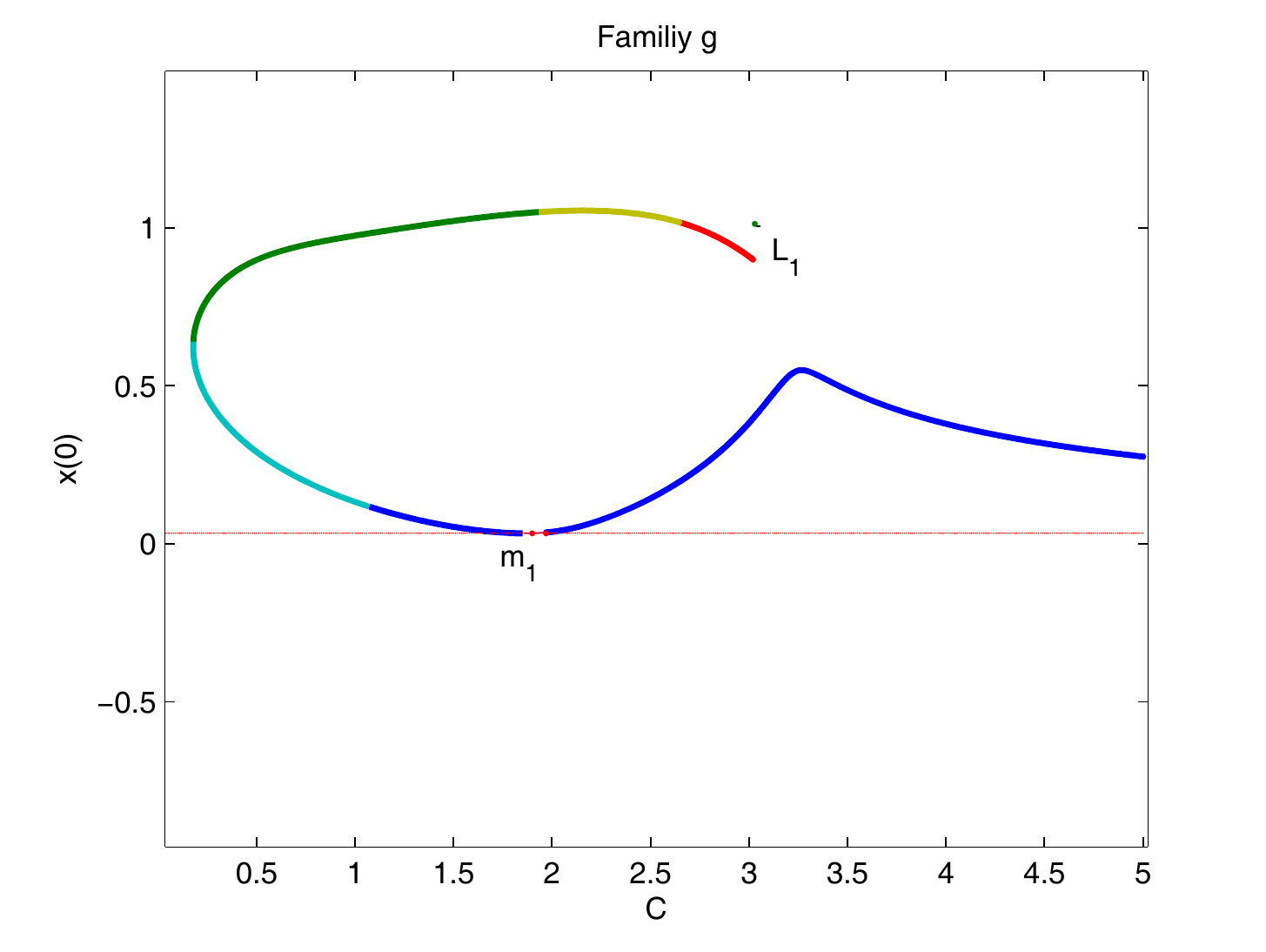}
\end{center}
\caption{Characteristic curve of family $g$.\label{familyg}}
\end{figure}
The first phase of this family stars with infinitesimal direct circular periodic orbits around $m_{1}$, the size of the orbits increases as the value of Jacobi constant $C$ decreases until a collision orbit is reached, the first phase of this family was established and studied in \cite{PapaI} for a little smaller value of $\mu$. Second phase starts with the retrograde orbit following the collision orbit which forms two loops as shown in Figure~ \ref{phasesg}. The inside loops now increase their size and the outside loops shrink as the value of $C$ deceases, both sets of loops become indistinguishable at the $fold$ point when $C\approx0.1797$. After this point the third phase starts, now the role of the loops is interchanged i,e., the inside loops shrinks and the outside loops expand as the value of $C$ increases.

Following the evolution of this phase we found that an orbit of collision with $m_{1}$ appears and the inside loops disappear, this is the beginning of fourth phase where the middle part of the orbits increases its size as the Jacobi constant increases, the termination of this phase (and of the whole family) are asymptotic orbits to $L_{2}$ ($L_{3}$ in \cite{PapaI}).
 More precisely, the value of $C$ oscillates in a small neighbourhood around the value of the Jacobi constant of the equilibrium point $L_{2}$ and the period tends to infinity as is predicted in theorem \ref{blue sky} similar to the Copenhagen category of R3BP.

\subsection{The family $f$ of retrograde orbits around $m_{1}$}

\begin{figure}[!h]
\begin{center}
\includegraphics[width=0.45\textwidth,height=0.45\textwidth]{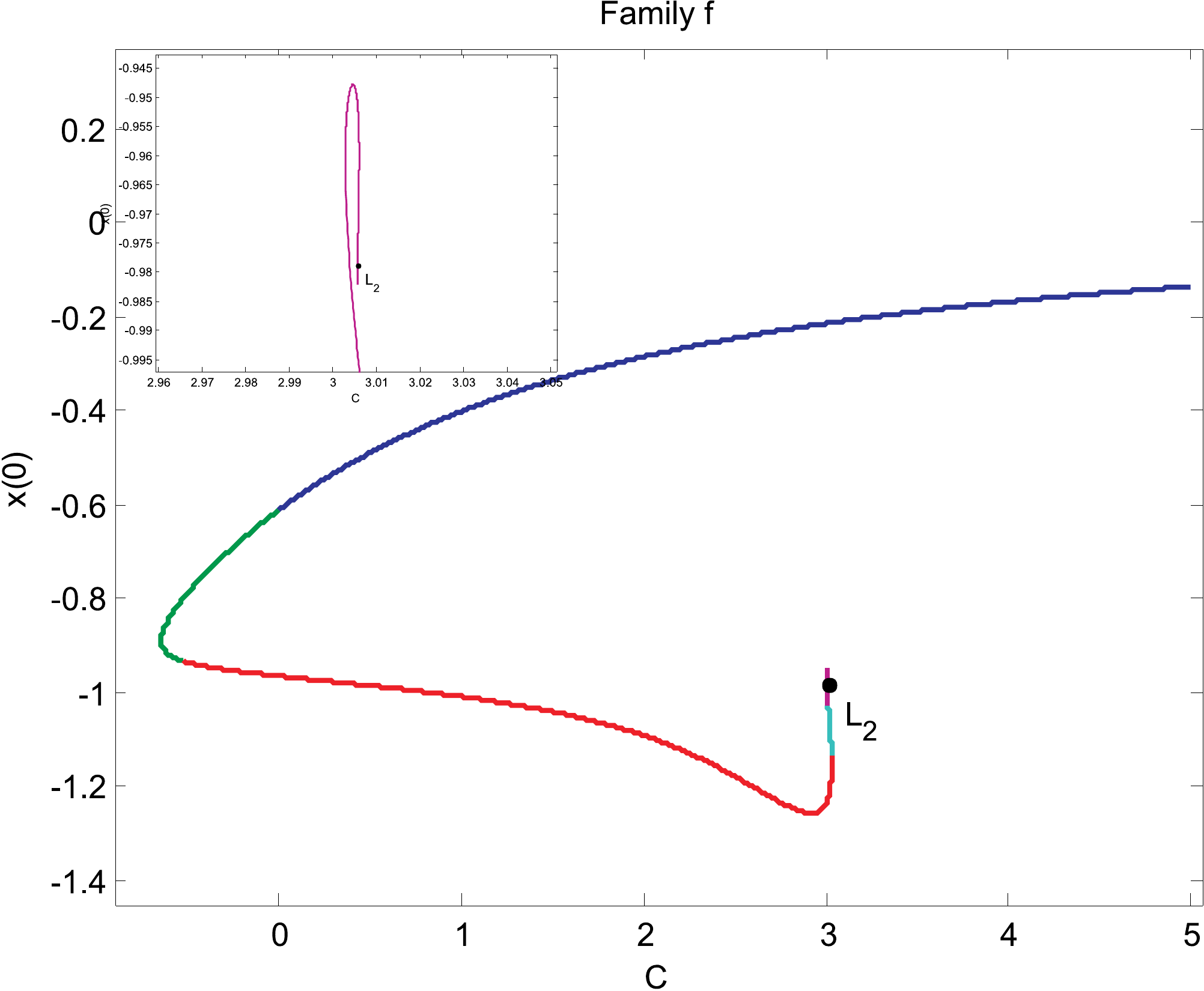}
\end{center}
\caption{Characteristic curve of family $f$, zoomed area indicates the end of family at the equilibrium point $L_{2}$.\label{familyf}}
\end{figure}
\begin{figure}
\begin{center}
\includegraphics[scale=0.35]{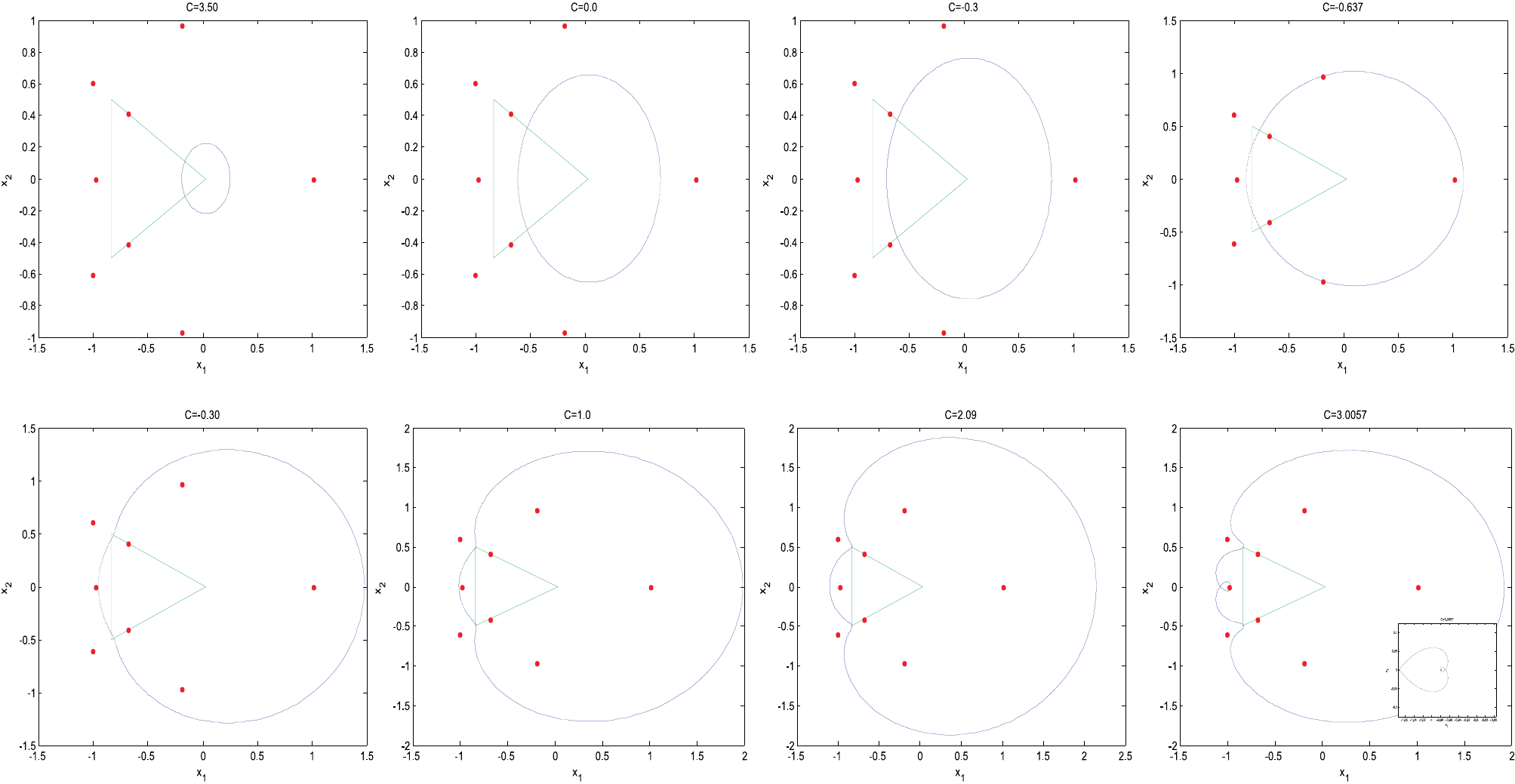}
\end{center}
\caption{Phases of family f.\label{phasesf}}
\end{figure}

The first phase of this family starts with infinitesimal retrograde circular periodic orbits around $m_{1}$, as in family $g$ the size of the orbits increase as the value of $C$ decreases monotonically until a $fold$ point is reached, this happens at $C\approx-0.6379$. This is the beginning of second phase, the periodic orbits still continue increasing their size but now these orbits tend to collision with the primaries $m_{2}$ and $m_{3}$, however this collision is never reached because the periodic orbits become asymptotic to $L_{2}$, as in family $g$.  See Figures~ \ref{familyf}, \ref{phasesf}.

\subsection{The family $a$ of retrograde orbits around $L_{1}$}
\begin{figure}[!h]
\begin{center}
\includegraphics[width=0.45\textwidth,height=0.45\textwidth]{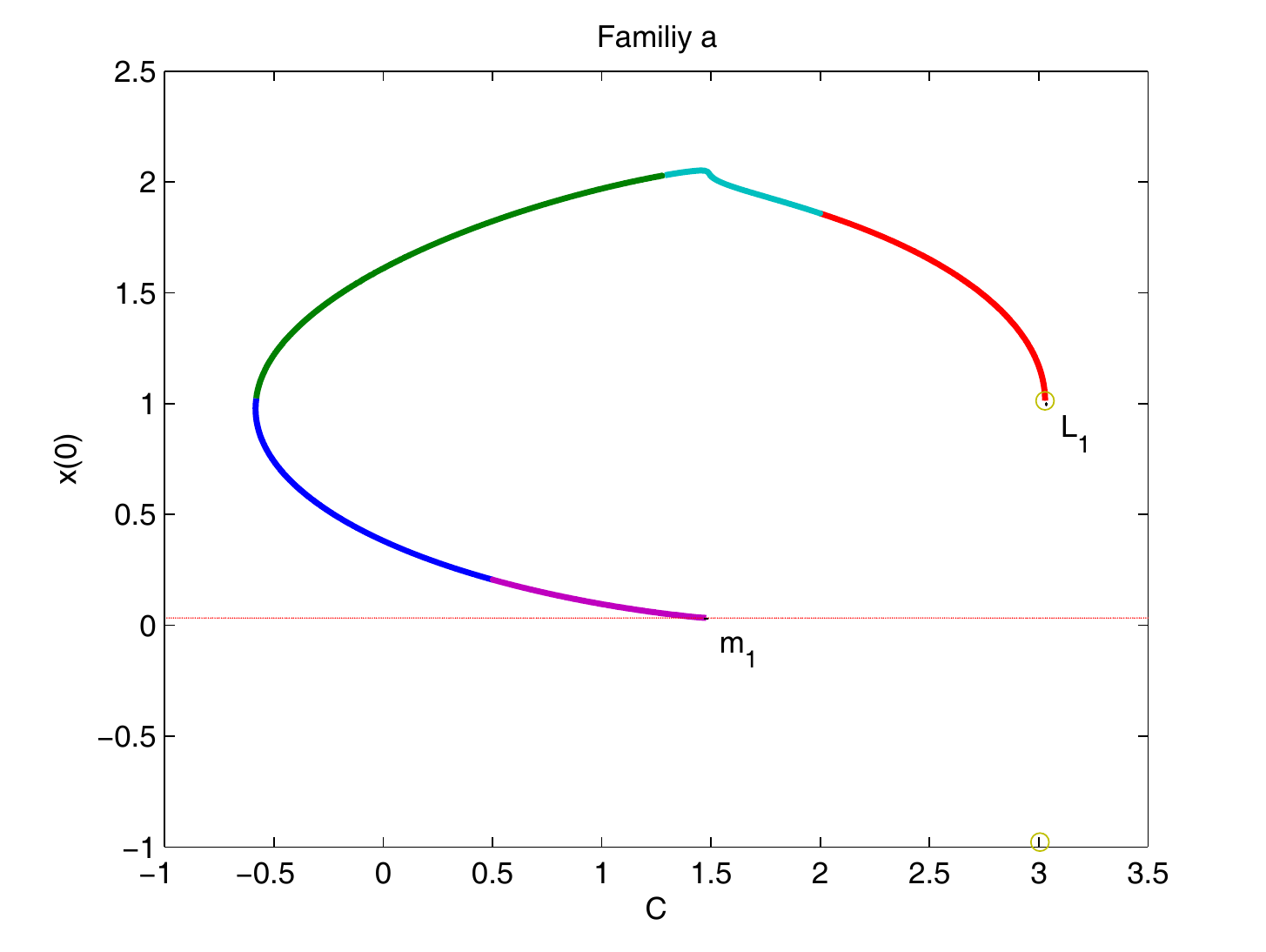}
\end{center}
\caption{Characteristic curve of family $a$.\label{characteristica}}
\end{figure}
\begin{figure}
\begin{center}
\includegraphics[scale=0.35]{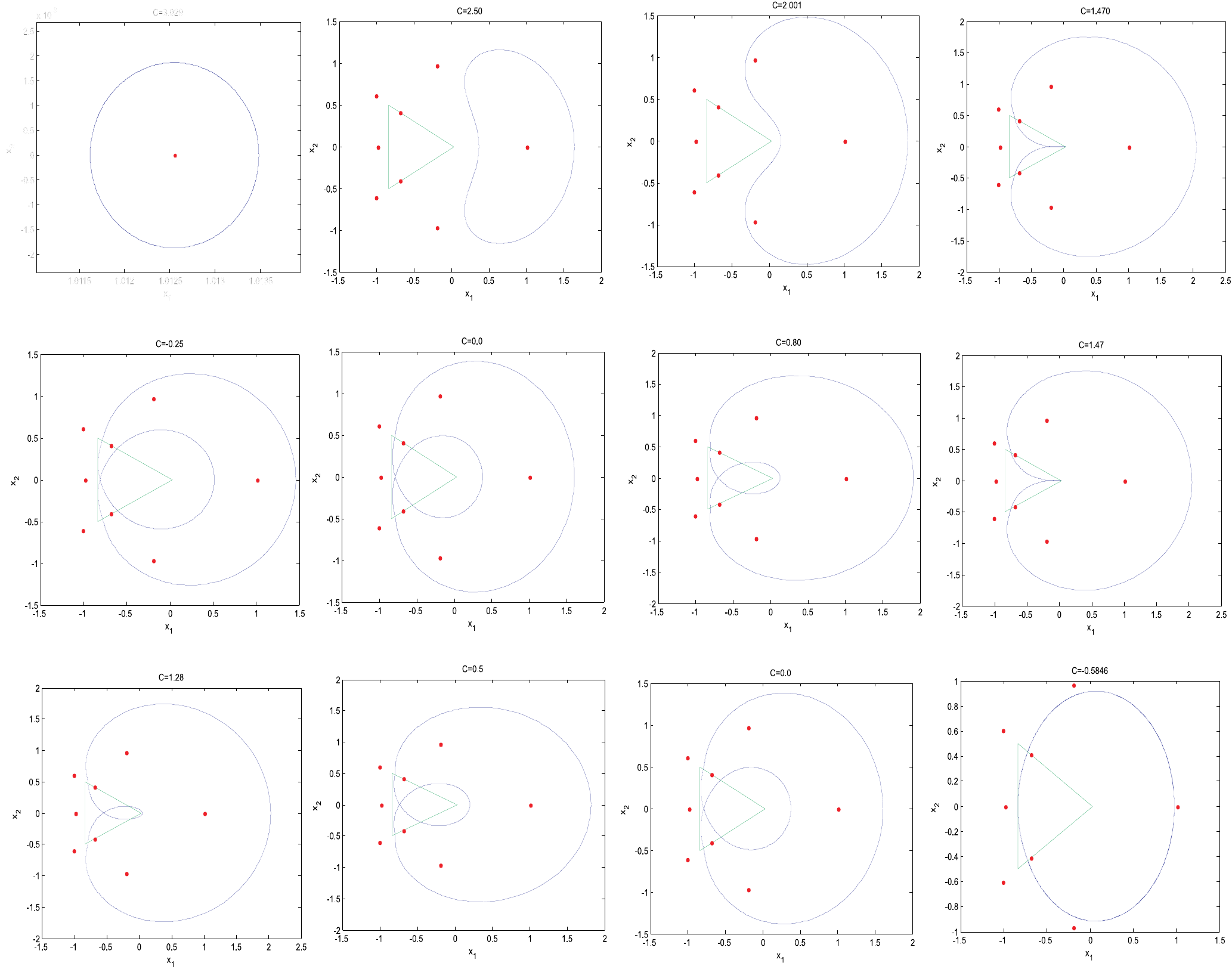}
\end{center}
\caption{Phases of family a.\label{phasesa}}
\end{figure}

The beginning of this family is provided by the Liapunov's center theorem, therefore in this family we continue retrograde periodic orbits around the equilibrium point $L_{1}$ ($L_{2}$ in \cite{PapaI})
for values of the Jacobi constant less than $C_{1}$. As the Jacobi constant decreases monotonically, the periodic orbits increase their size until a collision orbit with $m_{1}$ is reached; this is the end of first phase. While the value of $C$ continues decreasing, a second loop appears in the orbits, this loop increases its size along this second phase until a new $fold$ point is reached when $C\approx-0.5846$. At this point the inner and outer loops become indistinguishable as in the previous families, after this fold point, both loops are interchanged and the new inside loop shrinks as the value of $C$ increases until a collision orbit with $m_{1}$ finishes the third phase, see Figures~ \ref{characteristica}, \ref{phasesa}.

\subsection{The family $g_4$}
\begin{figure}[!h]
\begin{center}
\includegraphics[width=0.45\textwidth,height=0.45\textwidth]{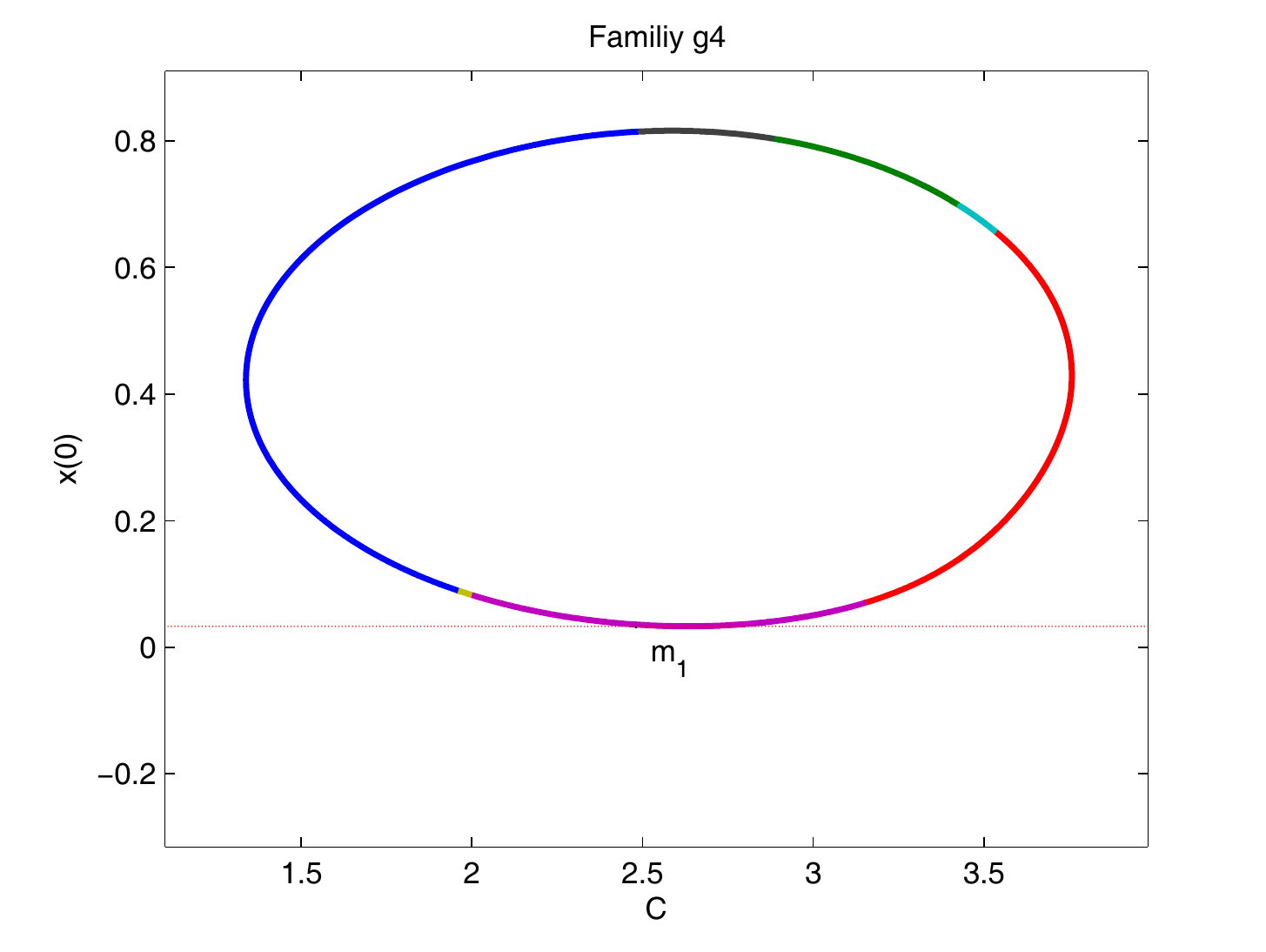}
\end{center}
\caption{Characteristic curve of the closed family $g_4$.\label{characterisitcg4}}
\end{figure}
\begin{figure}
\begin{center}
\includegraphics[scale=0.35]{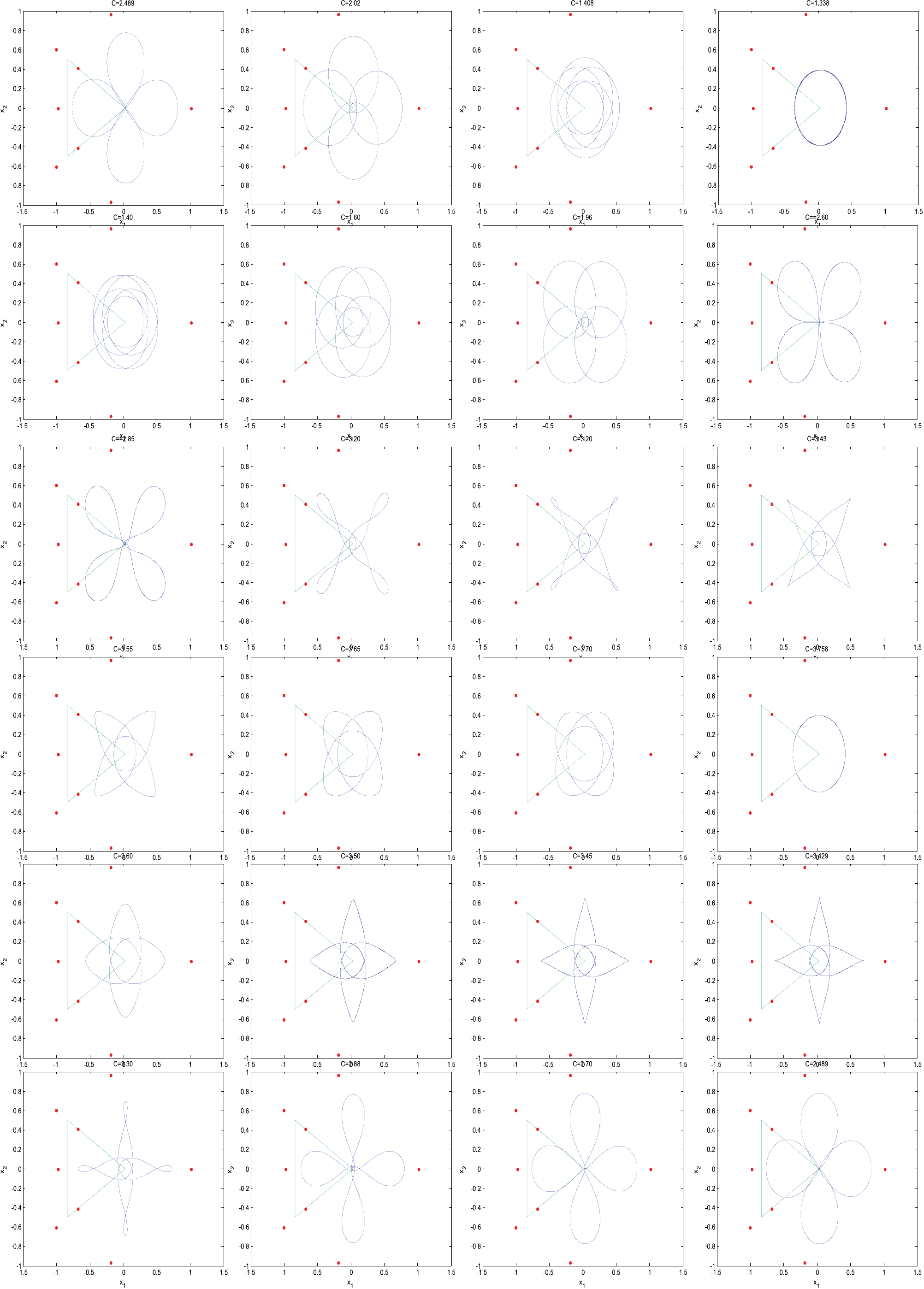}
\end{center}
\caption{Phases of family $g_{4}$.\label{phasesg4}}
\end{figure}

This family is the first example in this paper of a closed family of periodic orbits, see Figures~ \ref{characterisitcg4}, \ref{phasesg4}. The first phase of this family is composed of periodic orbits forming four loops around $m_{1}$, following the evolution of first phase we found that as $C$ decreases the four loops shrink until they become indistinguishable i,e; a $fold$ point is reached at $C\approx1.3381$, the second phase starts when these loops separate each other, as $C$ increases the periodic orbits tend to collision with $m_{1}$ i,e; the intersections between the loops tend to $m_{1}$. After this collision, the third phase starts. The periodic orbits change multiplicity because four inside loops around $m_{1}$ appear together with four outside loops. As the value of $C$ increases, the inside loops increase their size while the outside loops shrink at same time until they disappear. This is the end of third phase.

Fourth phase begins when the loops around $m_{1}$ increase their size as $C$ continues increasing monotonically. These loops become indistinguishable at the value $C\approx3.7581$ and a new $fold$ point is reached. After this $fold$ point the fifth phase begins, as expected inside and outside loops interchange and the peak of the outside loops become non-smooth i,e; no more orthogonal intersections with $x$-axis exist, this is the end of fifth phase. As $C$ decreases the loops of the orbits shrink to collision with $m_{1}$ and the sixth phase ends. We observe that the first phase of this family starts after this collision, therefore the family is closed.

\subsection{The family $g_6$}
\begin{figure}[!h]
\begin{center}
\includegraphics[width=0.45\textwidth,height=0.45\textwidth]{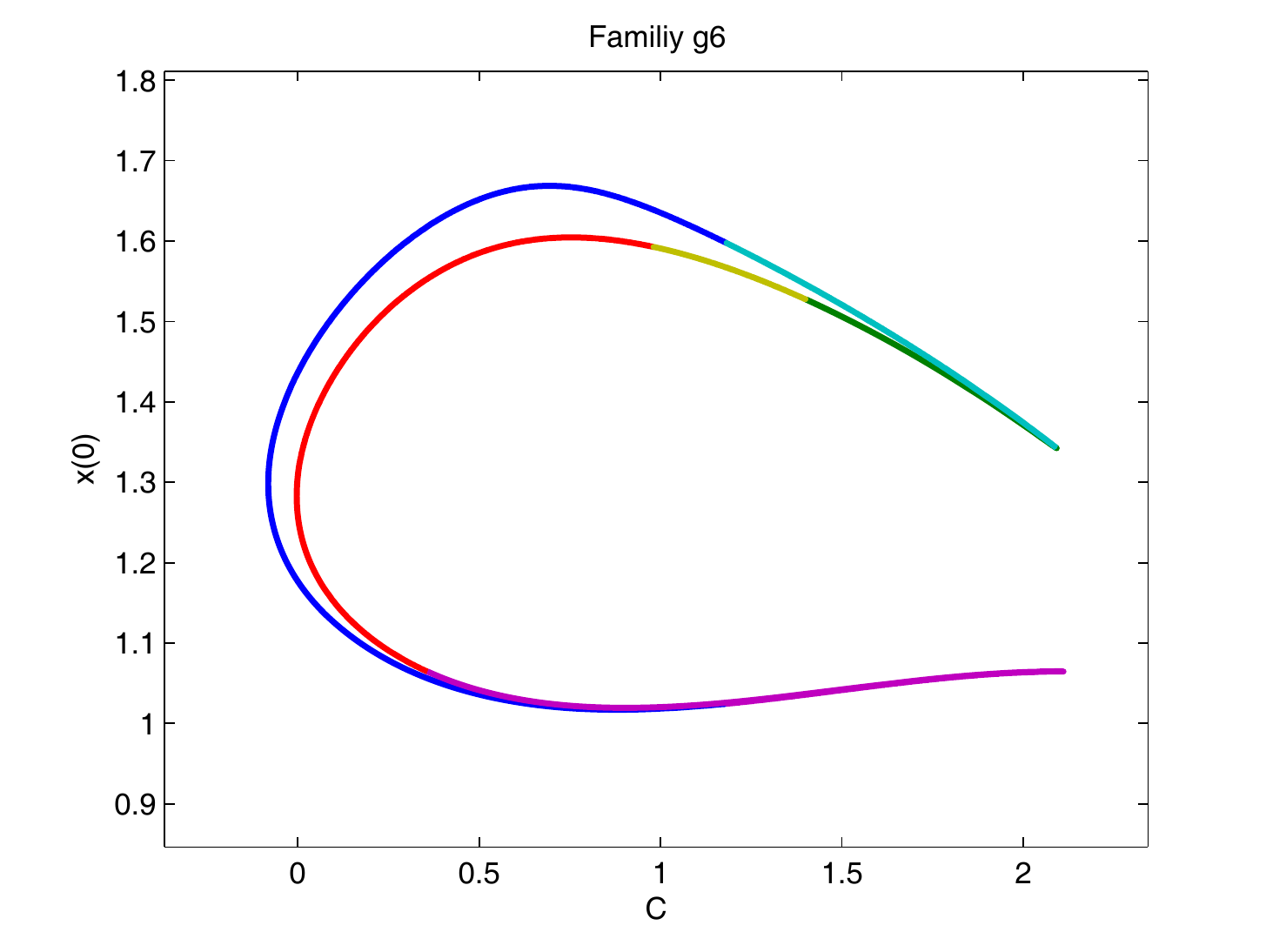}
\end{center}
\caption{Characteristic curve of the closed family $g_6$.\label{characterisitcg6}}
\end{figure}
\begin{figure}
\begin{center}
\includegraphics[scale=0.35]{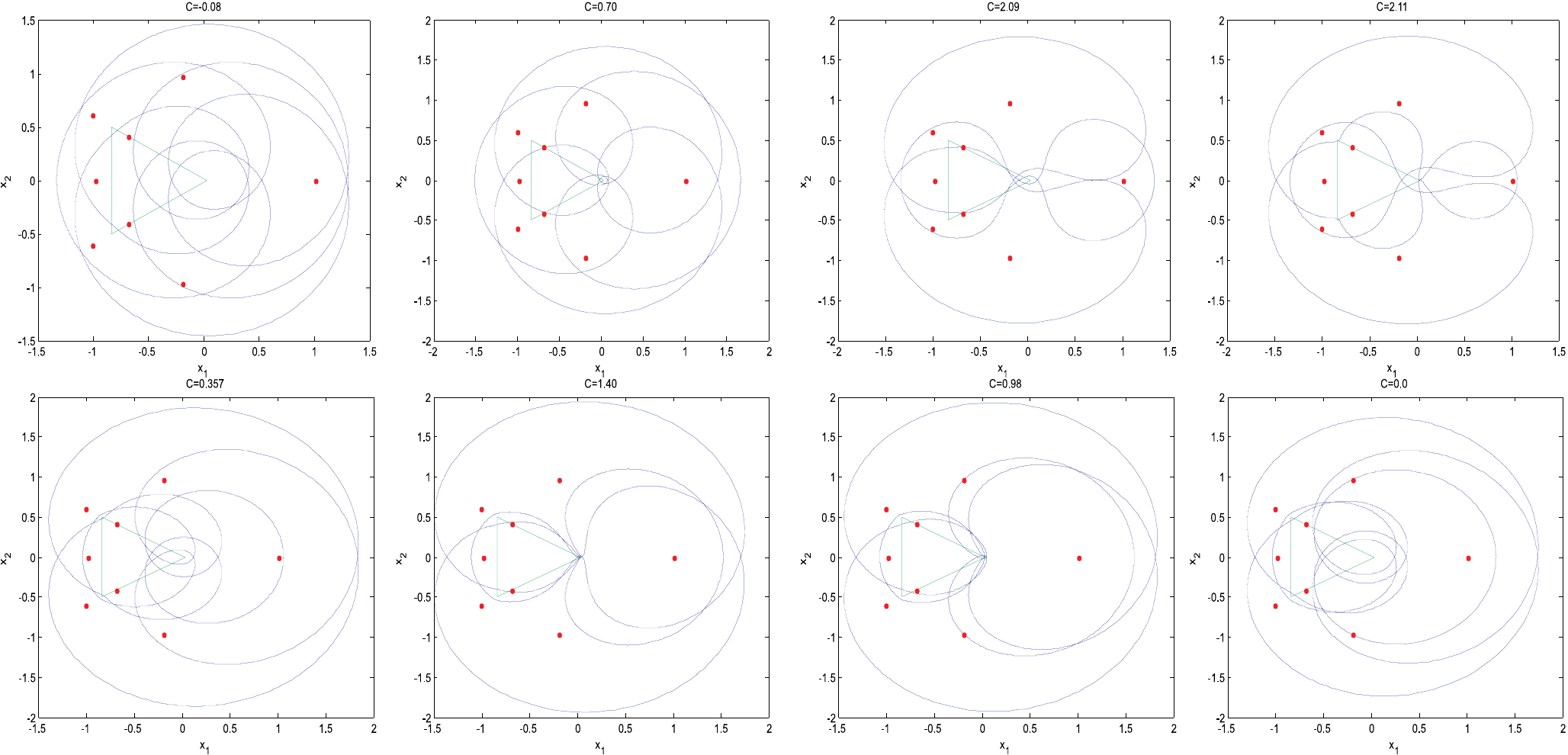}
\end{center}
\caption{Phases of family $g_{6}$.\label{phasesg6}}
\end{figure}

This family is another example of a closed family of periodic orbits. The behavior of its orbits is more complicated than that of previous families see Figure~ \ref{characterisitcg6}, \ref{phasesg6}. We have named $g_{6}$ this family because there exist a section of this family where the orbits form 6 loops around the primary $m_{1}$ but following the evolution of this family, we find that the orbits show a complicated behavior, therefore is not clear at all how to classify in phases this family, for simplicity we use the term ``first" phase to the outer section of the characteristic curve between the return points in Figure~ \ref{characterisitcg6}, and second phase to the inner section between the return points. In each case we show representative orbits in Figure~ \ref{phasesg6}.

\subsection{The family $m$ of retrograde orbits around $m_{1}$, $m_{2}$ and $m_{3}$}
\begin{figure}[!h]
\begin{center}
\includegraphics[width=0.45\textwidth,height=0.45\textwidth]{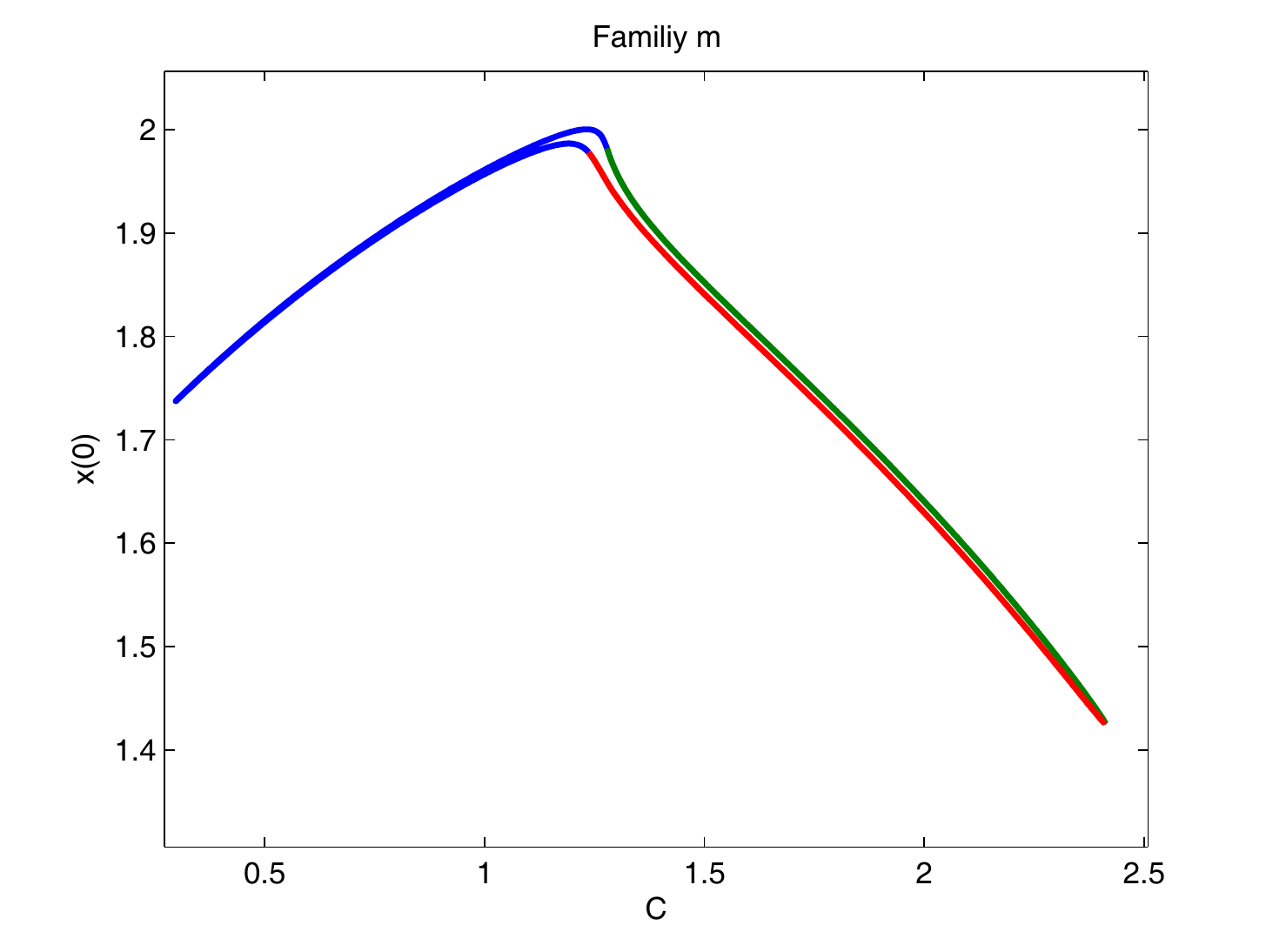}
\end{center}
\caption{Characteristic curve of the closed family $m$. \label{characteristicm}}
\end{figure}
\begin{figure}
\begin{center}
\includegraphics[scale=0.35]{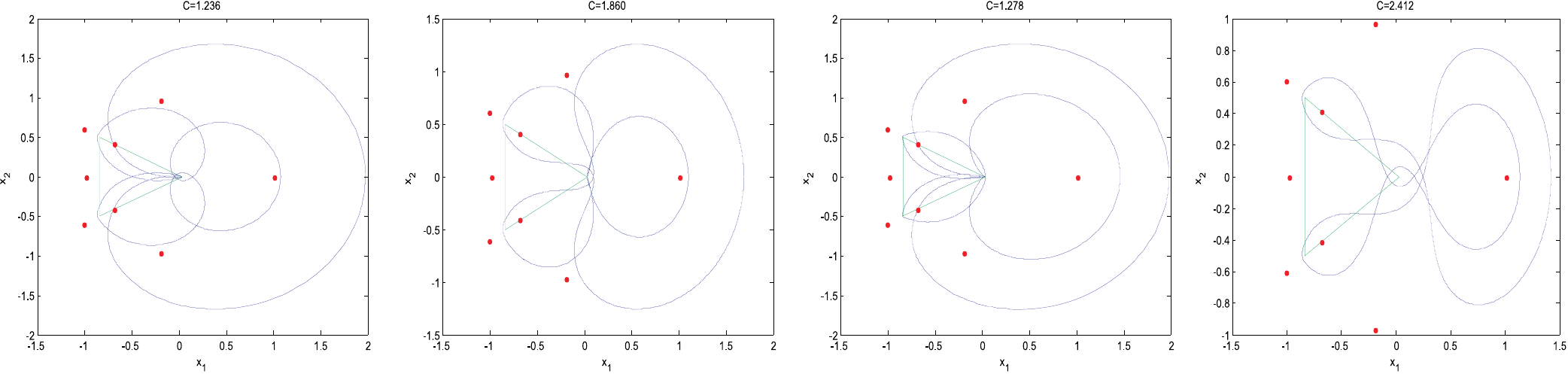}
\end{center}
\caption{Phases of family m.\label{phasesm}}
\end{figure}

Family $m$ consists of retrograde periodic orbits around the three primaries but these orbits do not surround the primaries simultaneously i,e., they form three loops, each loops surrounds one primary as can be seen in the Figure~ (\ref{phasesm}). The characteristic curve of the family is shown in Figure~ \ref{characteristicm} The three loops increase and decrease their size while the family is followed; however collision with the primaries never is reached although the orbits are close to collision when the loops decrease their size. This behavior is cyclic because the family is closed.

\subsection{The family $j$ of retrograde periodic orbits around $m_{2}$ and $m_{3}$}
\begin{figure}[!h]
\begin{center}
\includegraphics[width=0.45\textwidth,height=0.45\textwidth]{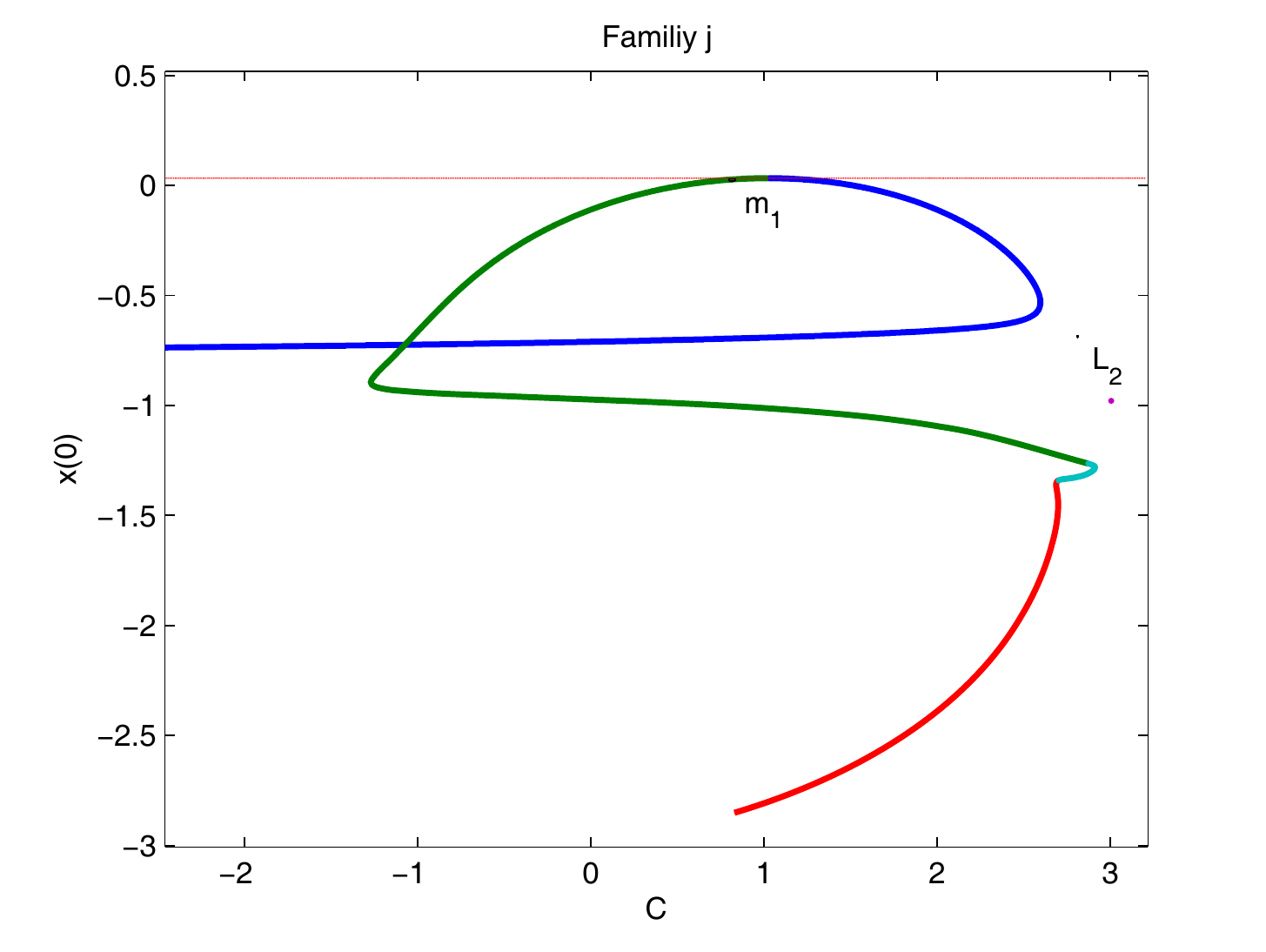}
\end{center}
\caption{Characteristic curve of the family $j$.\label{characteristicj}}
\end{figure}
\begin{figure}
\begin{center}
\includegraphics[scale=0.35]{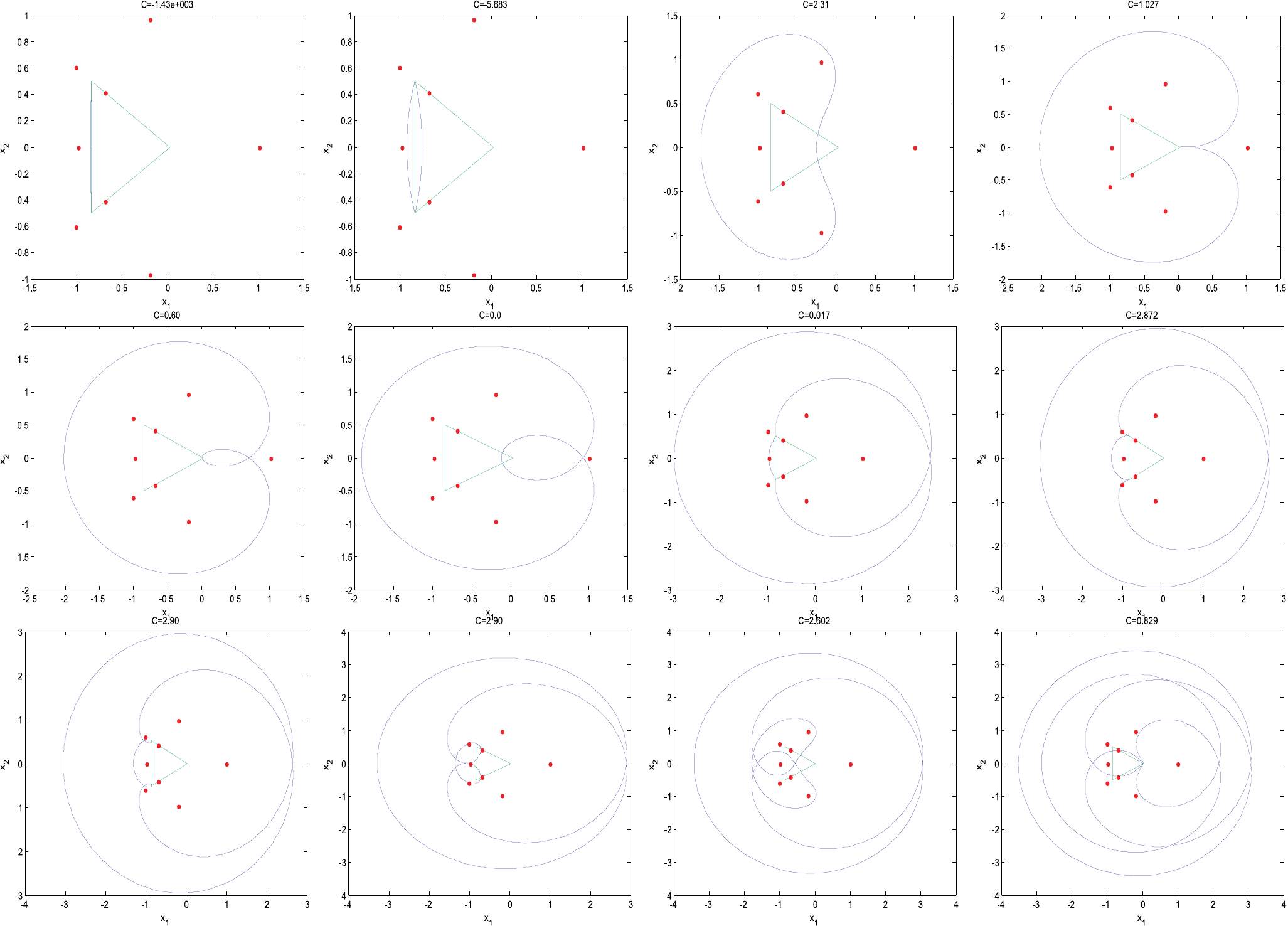}
\end{center}
\caption{Phases of family $j$.\label{phasesj}}
\end{figure}

The first phase of this family by one side tends to collision with the primaries $m_{2}$ and $m_{3}$, such collisions are reached for great values (negatives) of the Jacobi constant, see Figures~ \ref{characteristicj}, \ref{phasesj}. As $C$ increases, the orbits increase their size around $m_{2}$ and $m_{3}$, the orthogonal intersection to the right of both primaries $m_{2}$ and $m_{3}$ tends to the primary $m_{1}$ and a collision orbit appears,  this is the end of the first phase. After this collision an inner loop appears as expected and the orbits become direct around $m_{1}$, this is the second phase, the inside and the outside loops increase their size and the orbits tend to collision with both primaries $m_{2}$ and $m_{3}$ until this collision happens. The third phase starts after this collision, two loops around $m_{2}$ and $m_{3}$  respectively appear in the orbits and as $C$ decreases these loops these loops increase around the primaries, the behavior of the orbits complicate while the family is followed.

Finally a new collision orbit with $m_{1}$ appears. We have decided to terminate the family at this point, the complicated behavior of the orbits and the long time of integration of the regularized equations forced us to stop the continuation at this point.

\subsection{The family $r_{2}$ of asymptotic orbits to $L_{2}$}
\begin{figure}[!h]
\begin{center}
\includegraphics[width=0.45\textwidth,height=0.45\textwidth]{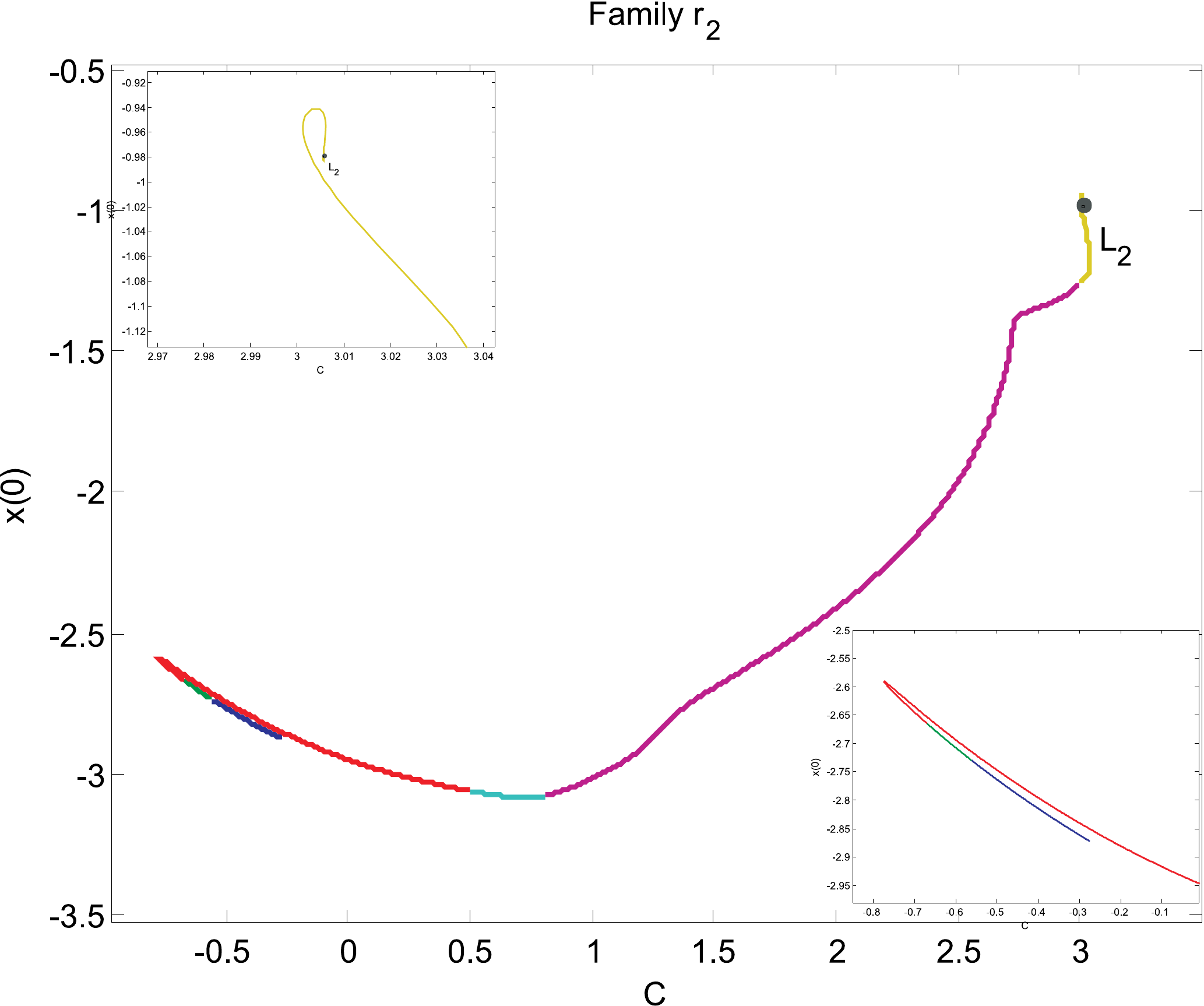}
\end{center}
\caption{Characteristic curve of family $r_{2}$, zoomed areas indicate the end of family at the equilibrium point $L_{2}$ and a return point of the family.\label{characteristicr2}}
\end{figure}
\begin{figure}
\begin{center}
\includegraphics[scale=0.35]{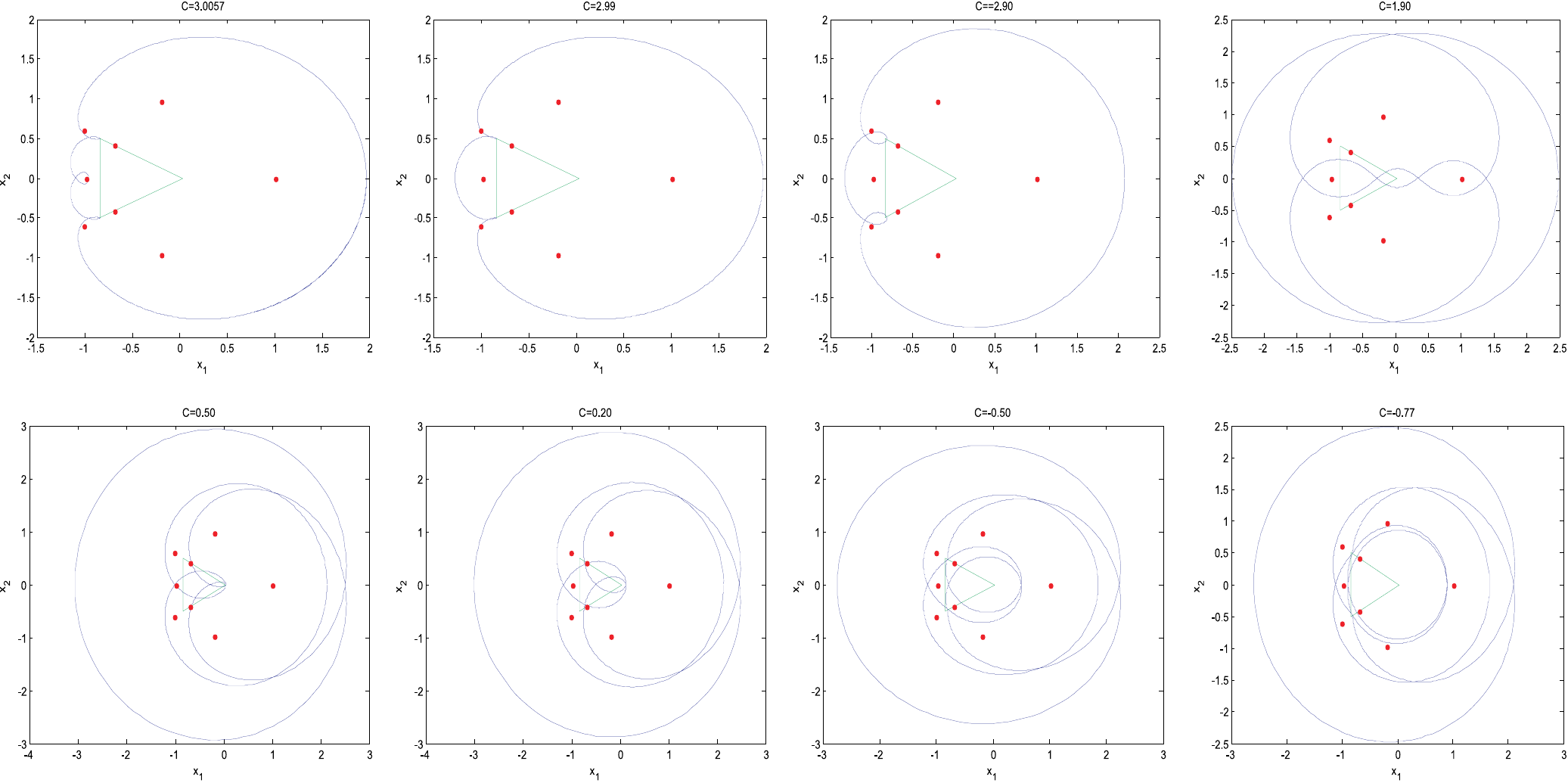}
\end{center}
\caption{Phases of family $r_{2}$.\label{phasesr2}}
\end{figure}

This family of periodic orbits has been named $r_{2}$ in analogy with the family $r$ of the Copenhagen category, the subscript indicates that the family is asymptotic to $L_{2}$ see Figures~ \ref{characteristicr2}, \ref{phasesr2}. The first phase of this family is composed by asymptotic periodic orbits to the equilibrium point $L_{2}$, these orbits form two loops surrounding the primaries $m_{2}$ and $m_{3}$ but while the periodic orbits go away from $L_{2}$ these loops shrink, therefore we have orbits close to collision now, however such collision never is reached. When the value of $C$ begins to decrease monotonically these loops increase its size and therefore the period of the orbits increase, it is interesting to note that the orbits become symmetric respect to $y$-axis (see Figure~ \ref{phasesr2}) however this symmetry disappears as $C$ continues decreasing, a collision orbits with $m_{1}$ terminates the first phase.

As expected, a new loop around $m_{1}$ appears in the orbits, here the second phase starts, following the evolution of the family we can see that the resulting loops of the collision increases its size as $C$ decreases monotonically and the orbits tend to collision with the primaries $m_{2}$ and $m_{3}$ . At $C\approx-0.7726$ we find a return point. We must emphasize that such collision is not reached although the mentioned loop continues increasing its size and therefore the orbits increase their period. This behavior continues as $C$ increases monotonically. As in family $j$ the long period of the orbits (long time of integration of the equations is needed) and the high instability of the orbits did not allow us to establish the end of this phase and therefore of the whole family.
\subsection{The family $j_{2}$ of retrograde periodic orbits around $m_{2}$ and $m_{3}$}
\begin{figure}[!h]
\begin{center}
\includegraphics[width=0.45\textwidth,height=0.45\textwidth]{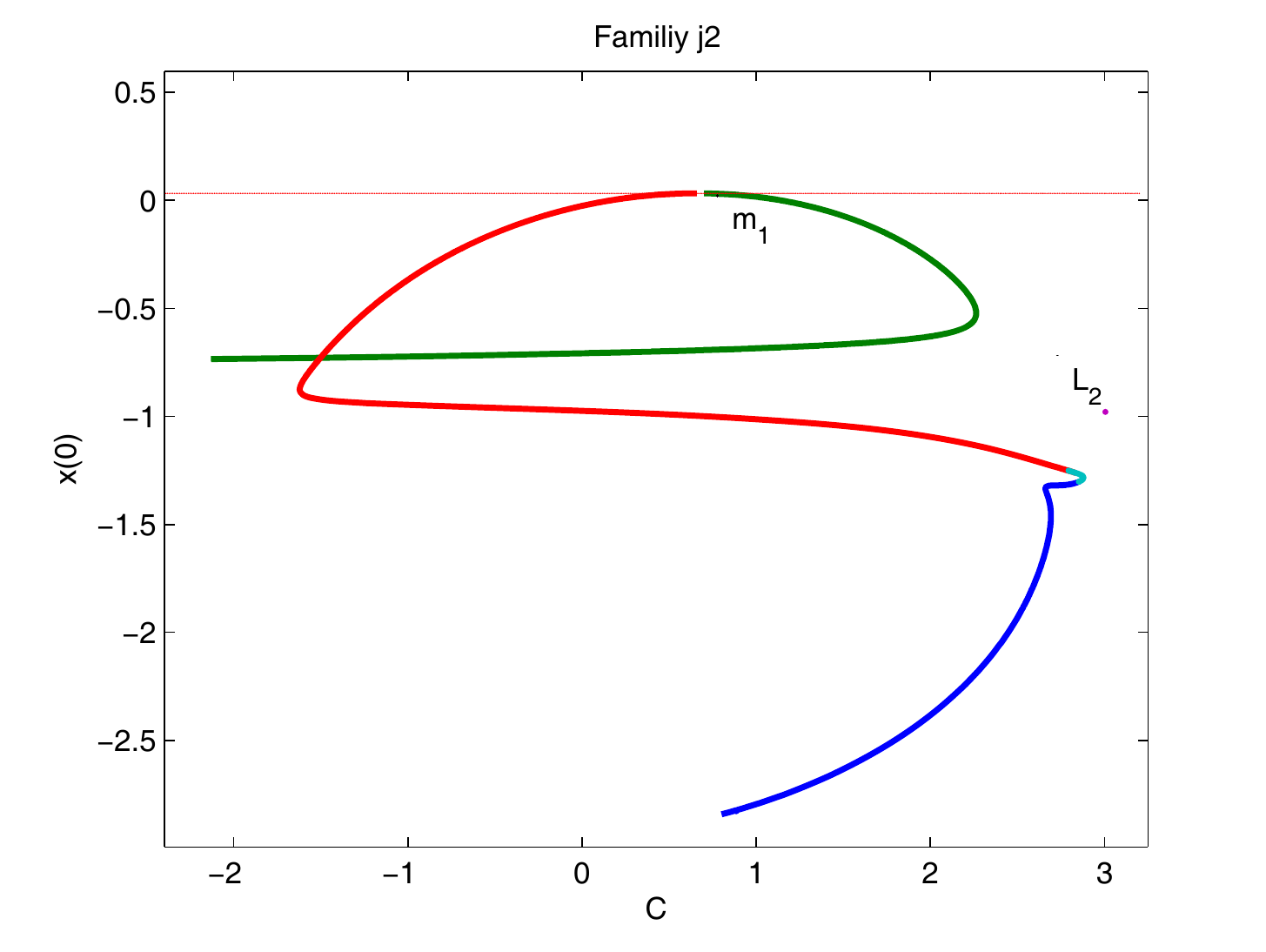}
\end{center}
\caption{Characteristic curve of the family $j_{2}$, note it is very similar to the one of the family j.\label{characteristicj2}}
\end{figure}
This last family named $j_{2}$ contains retrograde periodic orbits around $m_{2}$ and $m_{3}$ as in family $j$ but this time the orbits have an extra loop around $m_{2}$ and $m_{3}$ and this loop surrounds the primary $m_{1}$ too (see Figures~ \ref{characteristicj2}, \ref{phasesj2}). The evolution of this family is very similar to the evolution of family $j$, in fact, its characteristic curve has the same form that the one of family $j$. We show representative orbits of each phase in Figure~ \ref{phasesj2}.

\section{Stability of families and critical points}
\subsection{Critical points}
\begin{figure}[!h]
\begin{center}
\includegraphics[width=0.45\textwidth,height=0.45\textwidth]{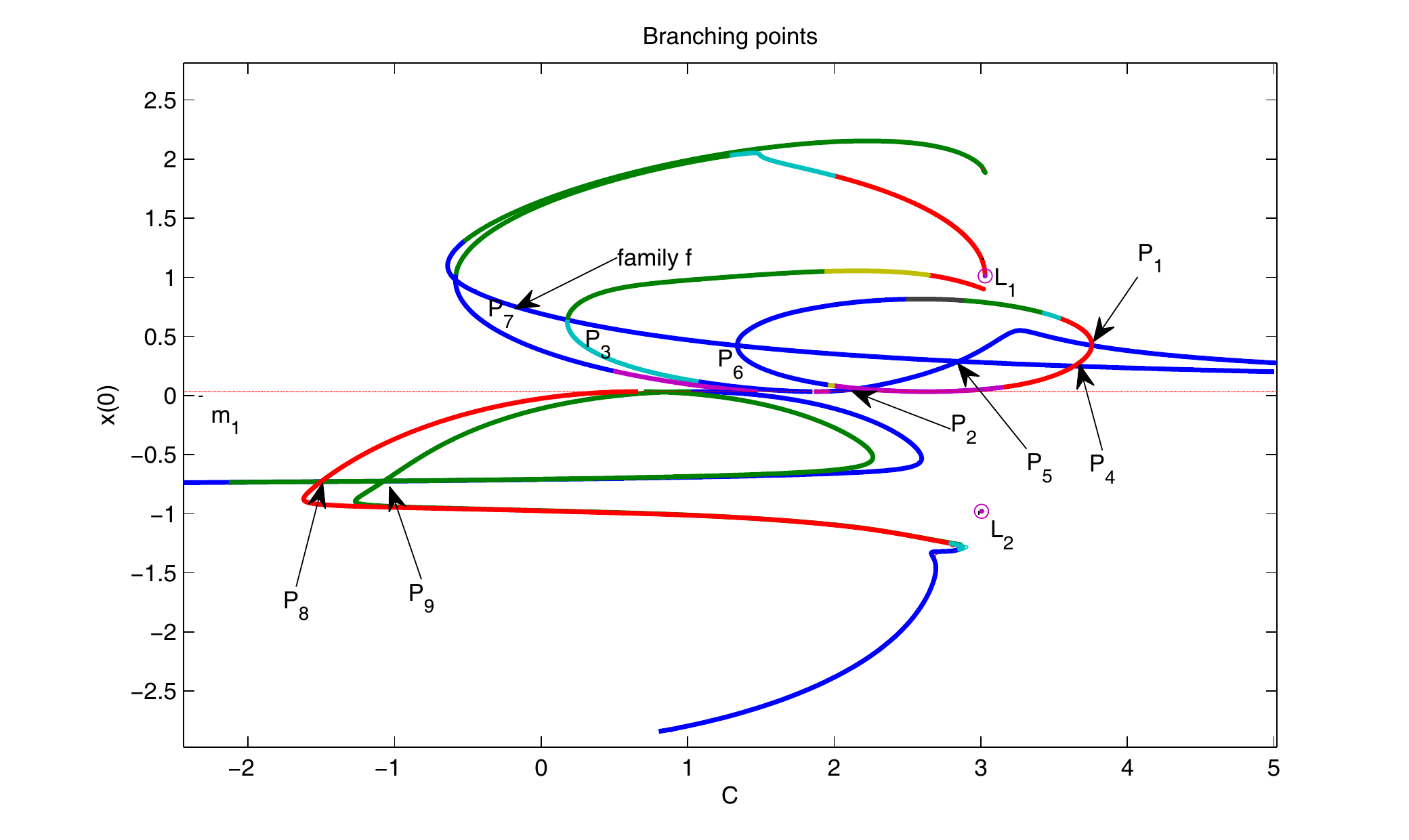}
\end{center}
\caption{Branching points between the families of periodic orbits\label{branchingpoints}}
\end{figure}
As can be seen in \cite{HenGen}, when we follow families of periodic orbits it can be found some especial points in such families called \textit{critical points} such as \textit{fold points} and \textit{branching points}, the last point is where two characteristic curves intersect, in 5.2 we will see that these points are in relation with stability changes, this facts do not ocurr by casuality, such phenomena was studied by example in \cite{HenII}. We must say that some families of periodic orbits were found through these points, in Figure~ \ref{branchingpoints} we show the branching points found in this study.

\begin{enumerate}
\item Family $g$ has 3 branching points; first one $P_{1}$ at $C\approx3.7581$ with family $g_{4}$, second one $P_{2}$ at $C\approx2.1662$ with family $g_{4}$ again, third one $P_{3}$ at $C\approx0.1797$ with family $f$. This last branching point happens at a fold point.
\item Family $f$ has 5 branching points, first one $P_{4}$ at $C\approx3.6364$ with family $g_{4}$, second one $P_{5}$ at $C\approx2.8481$ with family $g$, third one $P_{6}$ at $C\approx1.3381$ with family $g_{4}$, fourth one $P_{3}$ with family $g$ already mentioned, fifth one $P_{7}$ at $C\approx-0.5846$ with family $a$ ($fold$ point of family $a$)
\item Families $j$ and $j_{2}$ intersect at two points, first one $P_{8}$ at $C\approx-1.5015$ and second one $P_{9}$ at $C\approx-1.073$.

\end{enumerate}
\subsection{Stability of families}
\begin{figure}[!h]
\begin{center}
\includegraphics[width=0.45\textwidth,height=0.45\textwidth]{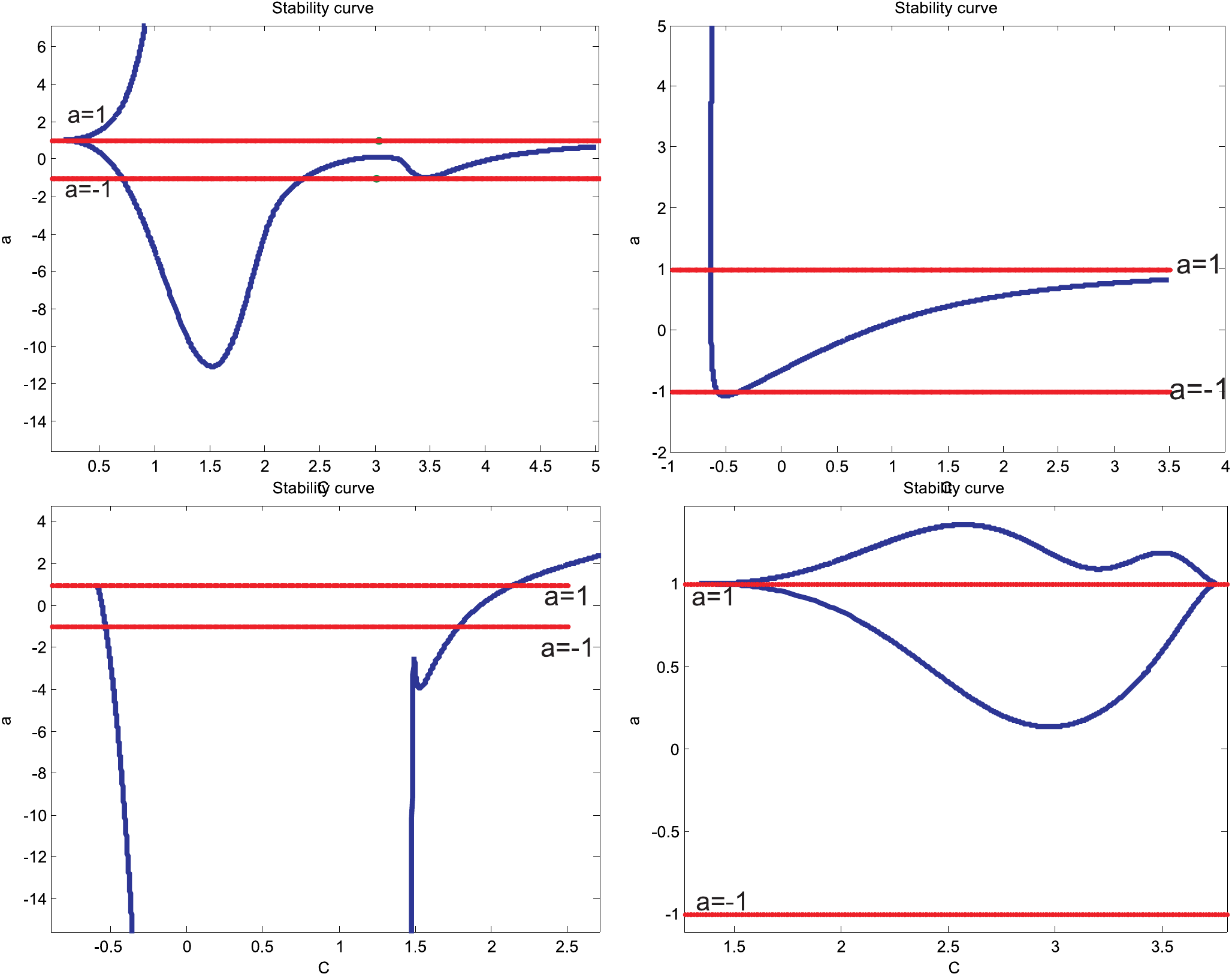}
\end{center}
\caption{Stability curve of family $g$ (up left), stability curve of family $f$ (up right), stability curve of family $a$ (down left), stability curve of family $g_{4}$ (down right).\label{stability1}}
\end{figure}

\begin{figure}[!h]
\begin{center}
\includegraphics[width=0.45\textwidth,height=0.45\textwidth]{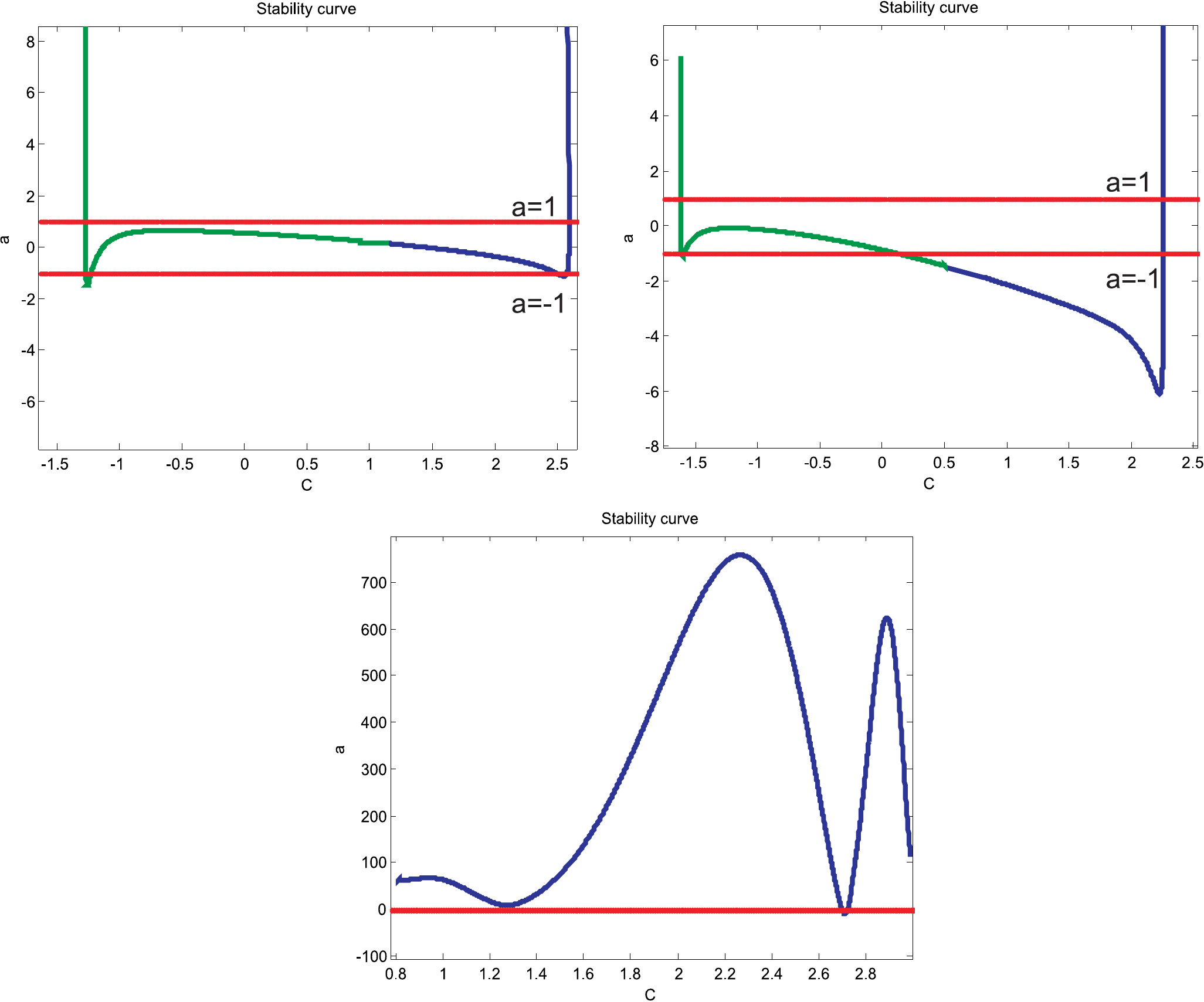}
\end{center}
\caption{Stability curve of family $j$ (up left), stability curve of family $j_{2}$ (up right), part of the end of the stability curve of family $r_{2}$ (center).\label{stability2}}
\end{figure}
In family $g$ we see that the first and second phase contain stable periodic orbits, at the end of second phase in the critical ($fold$) point we have that $\vert a\vert=1$ as expected, after this critical point all orbits become unstable, at the end of this family we observe strong oscillations of the sign of $a$ and therefore between the stable and unstable areas as is predicted by the ``blue sky catastrophe" termination. In family $f$ we have that the first phase of this family posses stable orbits as is shown in the Figure~ \ref{stability1}, at the $fold$ point we have $\vert a\vert=1$ again, after this point all orbits are unstable until the termination of family is reached and as in family $g$ strong oscillations between stable and unstable orbits are observed.

In family $a$ almost all orbits are unstable but we have 3 small regions where the orbits are stable, one of them is in the first phase, the second one is at the end of second phase; in the $fold$ point we have of course $\vert a\vert=1$. After this point the stability index starts to decrease and therefore we find a stability region at the beginning of third phase. The family $g_{4}$ is a closed family, a ``half" of its orbits are stable and a ``half" are unstable, this affirmation can be seen in Figure~ \ref{stability1}, at the two $fold$ points we have $\vert a\vert=1$ and between these two points the stability changes occur.

In family  $g_{6}$ almost all orbits are unstable, only a very little region around the two $fold$ points has stable orbits. Family $m$ presents a similar behavior as in family $g_{6}$. Orbits in family $r_{2}$ are strongly unstable except at the termination of family where oscillatory behavior between areas of stability and instability is observed. In the family $j$ we found a stability region between the values $C\approx-1.27$ and $C\approx2.59$, see Figure~ \ref{stability2} in this interval we have 3 critical points on the characteristic curve, and this can be seen on the stability curve where we have values for which $\vert a\vert=1$. For the family $j_{2}$ we have its stability curve behaves similar to the one of the family $j$, see Figure~ \ref{stability2}.

\section{Examples of periodic orbits}
\begin{figure}
\begin{center}
\includegraphics[scale=0.35]{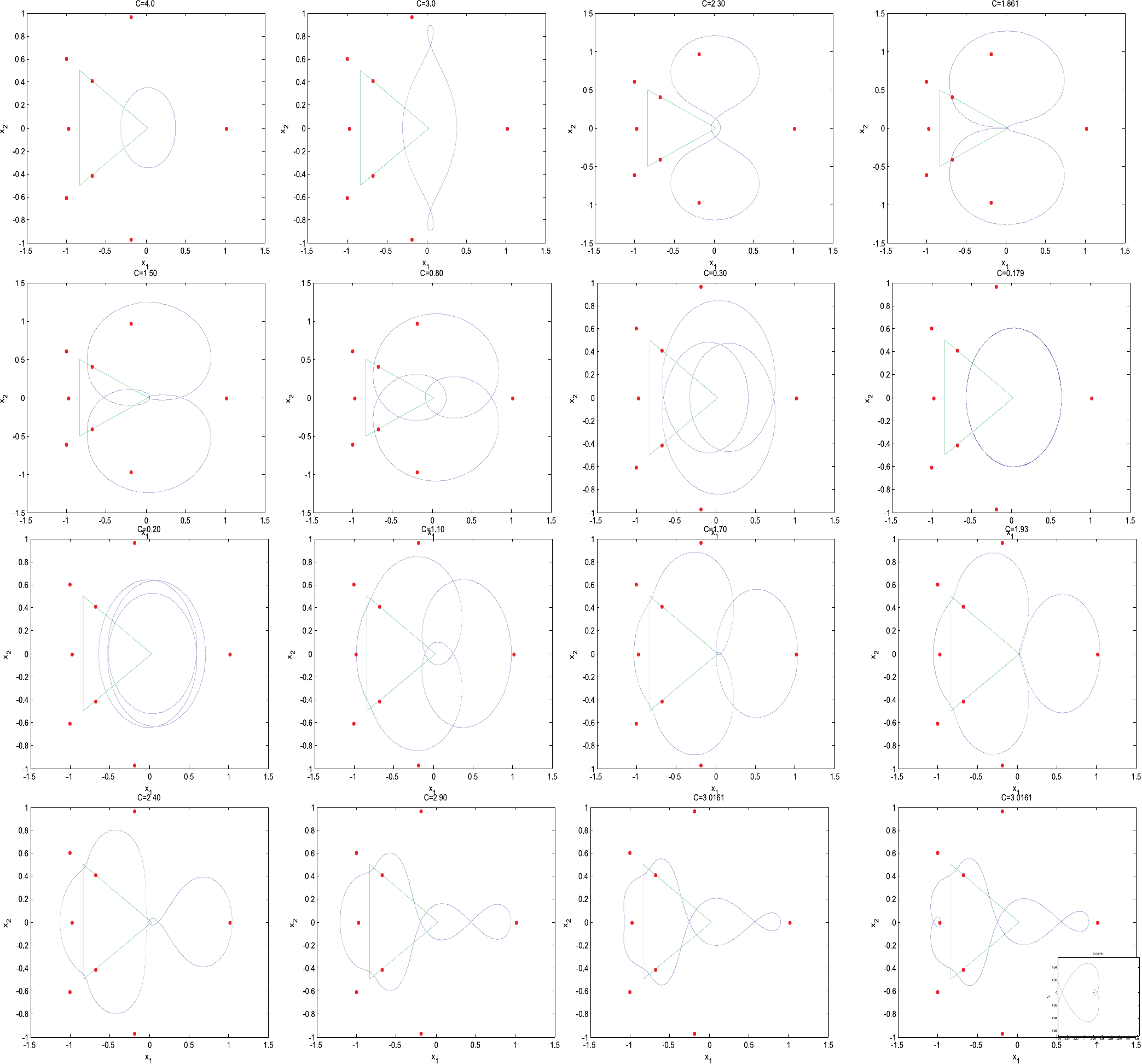}
\end{center}
\caption{Phases of family g.\label{phasesg}}
\end{figure}

\begin{figure}
\begin{center}
\includegraphics[scale=0.35]{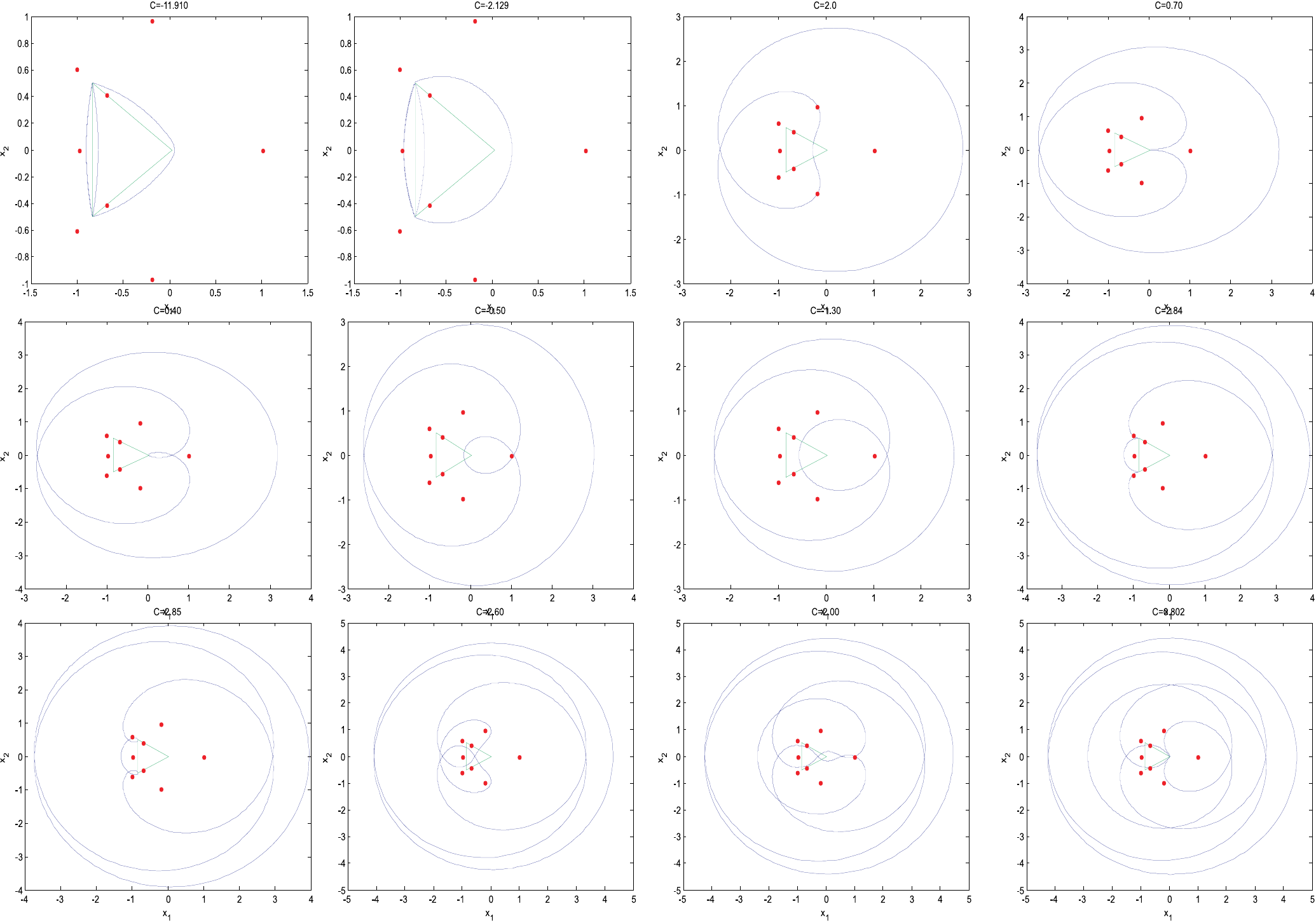}
\end{center}
\caption{Phases of family $j_{2}$.\label{phasesj2}}
\end{figure}

In Figures~ \ref{phasesg}, \ref{phasesj2} we show representative orbits of each phase of the nine found families of periodic orbits. Each row in the following figures represent a phase of the family, the first row represent the first phase of the family, second row represent the second phase etc.

\section{Conclusions and remarks}
In this work we extend the previous study of \cite{PapaI} and we include $5$ new families of periodic orbits. We performed the numerical continuation beyond different collisions with the primaries and show how some of them end up in a homoclinic orbit as predicted by the ``blue sky catastrophe" termination principle. In the former section we show representative orbits for each phase of the families, the orbits are mainly retrograde but there are families consisting of direct periodic orbits. We have also studied the stability of each family of periodic orbits and we show the characteristic curve of seven families where a considerable number of stable periodic orbits were found. Our main results can be summarized as follows:
\begin{enumerate}
\item The families $g$, $f$, $a$, $g_{4}$, $j$, $j_{2}$, and $r_{2}$ present a large number of stable periodic orbits, all the orbits in the families $g_{6}$ and $m$ are unstable except in very small regions.
\item The families $g$, $f$ and $r_{2}$ are asymptotic to the equilibrium point $L_{2}$, i,e; they end according to the ``Blue Sky Catastrophe" termination principle.
\item The families $g_{4}$, $g_{6}$ and $m$ are closed families.
\item The families $g$, $f$, $a$, $g_{4}$, $j$, and $j_{2}$ present \textit{Branching points} between their characteristic curves.
\item Almost all the families consist of retrograde periodic orbits, except families $g$, $g_{4}$ where have direct ones.
\end{enumerate}
\textbf{Acknowledgements} Author Burgos--Garc\'ia has been supported by a CONACYT fellowship of doctoral studies.
\newpage
\bibliographystyle{spr-mp-nameyear-cnd}

\bibliography{biblio-u1}

\end{document}